\documentclass[reqno,12pt] {amsart}

\newtheorem{theorem}{Theorem}[section]
\newtheorem{lemma}[theorem]{Lemma}

\theoremstyle{definition}
\newtheorem{definition}[theorem]{Definition}
\newtheorem{assumption}{Assumption}[section]

\theoremstyle{remark}
\newtheorem{remark}[theorem]{Remark}

\newcommand{\mysection}[1]{\section{#1}
\setcounter{equation}{0}}

\newcommand{\bR}{\mathbb R}

\newcommand\cL{\mathcal{L}}
\newcommand\cM{\mathcal{M}}

\newcommand{\sign}{\text{\rm sign}\,}

\begin{document}
\title[A priori estimates] {A priori estimates
of smoothness of solutions to  difference  Bellman equations 
  with linear and quasilinear operators}

\author[N.V. Krylov]{N.V. Krylov}

\address
{127 Vincent Hall, University of Minnesota, Minneapolis, MN
55455, USA}
\thanks{The work  was partially supported by
NSF Grant DMS-0140405}
\email{krylov@math.umn.edu}



\begin{abstract}
A priori estimates for finite-difference
approximations for the first and second order
derivatives are obtained for solutions
of parabolic equations described in the title.
\end{abstract}

\maketitle

\mysection{Introduction}

The goal of this article is to prove
a priori estimates for solutions
of finite-difference approximations
of parabolic Bellman equations with linear
and quasilinear operators.
In the latter case the nonlinear operator defining the equation
is still supposed to be convex with respect to the second-order
derivatives of the unknown function.
We present estimates for the
finite-difference approximations of the first
and second order spatial derivatives.
In particular, our results cover
finite-difference approximations for degenerate
quasilinear parabolic equations. As far as we are aware
these are the first results for such equations.
The main parts of the linear and quasilinear
operators entering Bellman equations
are assumed to be {\em linear\/}  
$a_{k}\Delta_{k}$ operators, that is written as 
a linear combination of pure second order derivatives
in certain directions that are common to all operators.
This assumption is always satisfied if the
equation is uniformly nondegenerate and is
generally  necessary
if we want to restrict ourselves to monotone
difference approximations and meshes
that are obtained from a fixed one by scaling
(see more about it in Remark \ref{remark 10.19.1} below).
Our results are valid for usual Bellman equations
and also for optimal stopping and impulse control problems
associated with them.

The motivation to obtain a priori estimates is the following.
There is an approach suggested in \cite{Kr97},
\cite{Kr99}, and \cite{Kr00}
to establishing the rate of convergence of $u_{h}$
to $u$ as $h\downarrow0$, where $u$ is the true solution,
 $u_{h}$ the solution
of finite-difference approximation of the same equation, and
$h$ typically is the mesh size. 
Two main ideas of this approach are
that the original equation and its finite-difference
approximation should play symmetric roles
and that one can ``shake the coefficients"
of the equation in order to be able to
mollify under the sign of nonlinear operator.

For elliptic Bellman equations with constant coefficients
and Lipschitz free terms
the first idea led to the rate of convergence of order
 $h^{1/3}$, 
for generic finite-difference approximations
and $h^{1/2}$ in the case of $a_{k}\Delta_{k}$ operators (see Remark 1.4
and Theorem 5.1 in \cite{Kr97}, also see
\cite{BJ1}). 
In contrast with the popular belief that
assuming more smoothness of the data does not lead 
to better rates of convergence, it is proved 
in \cite{DK1} that if the free terms are in $C^{1,1}$,
then the rate is at least
$h$ for constant coefficient $a_{k}\Delta_{k}$ case and $h^{2}$
for equations with better structure.

The second idea was introduced to treat equations with
variable coefficients and led to quite
satisfactory  error bounds for $u-u_{h}$ from an ``easy" side
(depending on how the equation is written
this can be either upper or lower estimate of $u-u_{h}$).
To get an estimate from the other side
on the basis of the idea of symmetry between
the approximating and the original equations
one needed to solve the following problem:

 (P) in the case of variable coefficients estimate
 how much the solution
of the finite-difference equation loses
 in the process
of shaking the equation. 

In the absence of solution
of the problem (P) the idea of symmetry 
was still useful but only
in obtaining some intermediate
estimates  (see, for instance,  \cite{BJ1} and   \cite{Kr99}) and
various approaches to getting the  
error bounds from the ``hard" side were developed. In addition to
the above cited papers the interested reader should
consult  \cite{BJ2}, \cite{BJ3}, and the references
therein. Note that for generic finite-difference approximations,
under the assumptions of Theorem 5.3 of \cite{Kr99} 
the result of \cite{BJ3} is the same $h^{1/3}$,
but the result of Theorem 5.4 of \cite{Kr99}
is improved from $h^{1/21}$ to $h^{1/7}$. The issue of 
solving the  problem (P) 
for generic finite-difference
approximation remains unsettled and 
it is not clear how
far off  $h^{1/7}$  
is  from the true rate.

The problem (P) was recently reduced to the problem
of estimating
the modulus of continuity of approximate 
solutions and solved in \cite{Kr06}
for $a_{k}\Delta_{k}$ case in which a sharp 
error bound of order $h^{1/2}$ was obtained.
The idea of symmetry worked again as in
the constant coefficients case.
This activity was continued in \cite{DK3},
where for the first time equations in domains
were treated,  and in \cite{DK2}, where
under various smoothness assumption
the rates $h^{1/2}$, $h$, and $h^{2}$
are obtained for {\em linear\/} degenerate
equations of $a_{k}\Delta_{k}$ form. For linear case
the rate $h^{1/2}$ was earlier obtained in \cite{Ja}
by a method close to a method from \cite{Kr00}
(Lemma 5.1 of \cite{Ja} is a version of
Theorem 2.1 of \cite{Kr00}).
However, this method does not allow 
one to get rates $h$ and $h^{2}$.

The main technical result of \cite{Kr06}
is the a priori estimate of the
derivative of $u_{h}$ with respect to $x$
stated as Theorem 5.2 and 
 proved by quite subtle estimates.
It turns out that there is a much easier method to prove 
 Theorem 5.2 of \cite{Kr06} which in addition 
carries over to much more
general equations with quasilinear operators
and to obtaining estimates for the second-order
finite differences of $u_{h}$. The method is almost as simple
as the one used in \cite{DK2} for
linear equations.

We present this new method here and concentrate only
on a priori estimates to keep the article within
reasonable limits. Once the a priori
estimates are obtained, one can follow familiar patterns
to get error bounds in various cases
of linear or quasilinear operators, degenerate
or weakly nondegenerate or else uniformly
nondegenerate,  with $C^{1}$
or $C^{1,1}$ coefficients. 
In particular, we hope to obtain first 
estimates on the rate of convergence 
in the case of degenerate quasilinear operators with Lipschitz
coefficients.
Our preliminary 
computations also show that  under the assumptions of Theorem
\ref{theorem 10.4.1}   the estimate
$|u-u_{h}|\leq N h^{2/3} $ holds in the elliptic case.
These and some other indicated below 
possible applications of our results
we intend to develop in the future.

Hongjie Dong and the referees of the paper
 made valuable comments on the first version
of it for which the author is sincerely
grateful.

\mysection{Setting and main results}
                                            \label{section 5.9.2}

Our first few results concern equations of the type
\begin{equation}
                                                    \label{10.26.1}
 F(\delta^{T}_{\tau}u,\Delta_{h,\ell_{k}}u,\delta_{h,\ell_{k}}u,u)=0,
\end{equation}
where
$$
F(\phi,q_{k},p_{k},\psi)=F(\phi, q ,p,\psi,t,x)
$$
\begin{equation}
                                                    \label{11.9.1}
=
\sup_{\alpha\in A}[r^{\alpha}(t)\phi+
\sum_{|k|=1}^{d_{1}}(a_{k}^{\alpha}(t,x)q_{k}
+b_{k}^{\alpha}(t,x)p_{k})-c^{\alpha}(t,x)\psi+f^{\alpha}(p,\psi,t,x)] ,
\end{equation}
$\delta^{T}_{\tau}u,\Delta_{h,\ell_{k}}u,\delta_{h,\ell_{k}}u$
are finite-difference approximations of the time derivative,
the pure second-order derivative in direction $\ell_{k}$,
and the first-order derivative in direction $\ell_{k}$,
respectively. Detailed description of the above objects now follow.

Let $A$ be a separable metric space, $d,d_{1}\geq1$  integers, and let
$$ 
r^{\alpha}=r^{\alpha}(t),\quad
a_{k}^{\alpha}=a^{\alpha}_{k}(t,x),\quad
 b_{k}^{\alpha}=b^{\alpha}_{k}(t,x),\quad
c^{\alpha}=c^{\alpha} (t,x),\quad
$$
 be real-valued {\em bounded\/}
functions of $(\alpha,t,x)$ 
defined on $A\times\bR\times\bR^{d}$ for $k=\pm1,...,\pm d_{1}$.
Also let some vectors $\ell_{k}\in\bR^{d}$ be defined
for $k=\pm1,...,\pm d_{1}$ and let 
\begin{equation}
                                                          \label{5.13.1}
T,h_{0} \in(0,\infty),\quad \delta\in(0,1],\quad
 K_{0},K_{1},K_{2},K_{3}\in[0,\infty),\quad m\in\bR
\end{equation}
 be some constants fixed throughout the article.
It is worth noting that   $\ell_{k}$, $k=\pm1,...,\pm d_{1}$,
are {\em not\/} supposed to form a basis in $\bR^{d}$
or even generate $\bR^{d}$.
This becomes crucial
when one proves the estimates of the first-order differences
of solutions with respect to parameters on which the coefficients
may depend. Notice also that the  lengths of $\ell_{k}$'s
can be different and some of them can be just zero
(and we will use this possibility later).
The constant $T$ gives us the time interval $[0,T)$, on which
the equation is investigated, $h_{0}$ ``calibrates"
the mesh-sizes in $x$ variable, the constant $\delta$ will appear
in various requirements of nondegeneracy.
The constant $K_{0}$ is the most basic one,
it is used in formulations of the very basic
assumptions. The constant $K_{1}$ is used to control
either the maximum magnitude of the solution
or its oscillation. The constant $K_{2}$ will appear
in our assumption on the growth of $f$ with respect to the ``gradient"
of the solution (see Assumption \ref{assumption 9.27.4} (ii),
which looks very much like the one
commonly used in the theory of quasilinear PDEs.
By the way, the author's efforts to use 
Assumption \ref{assumption 9.27.4} (iii), stated similarly,
  failed.)
The constant $K_{3}$ is used to control various quantities
having lesser impact on  our results than those
controlled by $K_{0},K_{1},K_{2}$.
Finally, the constant $m$ is used to extract various
results, which in the theory of parabolic PDEs
one gets after replacing $u(t,x)$ with $u(t,x)e^{mt}$.

For any vector $l\in\bR^{d}$, $\eta,\tau>0$, and function $u$ introduce
$$
\delta_{\eta,l}u(x)=\frac{u(x+\eta l)-u(x)}{\eta},\quad
\tau_{T}(t)=\tau\wedge(T-t)^{+},
$$
$$
\delta^{T}_{\tau}u(t,x)=\frac{u(t+\tau_{T}(t),x )-u(t,x)}{\tau},
\quad
\delta _{\tau}u(t,x)=\frac{u(t+\tau ,x )-u(t,x)}{\tau},
$$
$$
\Delta_{\eta,l}u(x)=\frac{u(x+\eta l)-2u(x)
+u(x-\eta l)}{\eta^{2}},
$$
where the notation $a^{\pm}=(1/2)(|a|\pm a)$ is used.
Observe that with the above definition of $\delta^{T}_{\tau}$
equation \eqref{10.26.1} makes perfect sense for $t<T$
for functions $u(t,x)$ defined only for $t\leq T$. We do not need
to extend $u$ beyond $T$ in order to compute
the finite-difference approximation of its derivative in time
for $t<T$.

\begin{assumption}
                                         \label{assumption 9.27.1}
(i) The functions $r^{\alpha}$, $a^{\alpha}_{k}$,
$b^{\alpha}_{k}$,
and $c^{\alpha} $ are continuous with respect to $\alpha$;

(ii) the functions  
$ b^{\alpha}_{k}$ satisfy the Lipschitz
condition with constant $K_{0}$ with respect to $x$;

(iii) the function $c^{\alpha} $ satisfies the Lipschitz
condition with constant
$K_{3}$ with respect to $x$;

(iv) we have
$$
\ell_{-k}=-\ell_{k},\quad a^{\alpha}_{-k}=a^{\alpha}_{k},\quad
|\ell_{k}|\leq K_{0} ,\quad  r^{\alpha}_{k}\geq0,\quad
a^{\alpha}_{k}\geq0 
$$
(for all values of the arguments and $k$).

\end{assumption}

An important feature of Assumption
\ref{assumption 9.27.1} is that no control
on the sizes of $r^{\alpha}$, $a^{\alpha}_{k}$,
$b^{\alpha}_{k}$,
and $c^{\alpha} $ is imposed (however, remember that
from the very beginning they are assumed to be bounded).

\begin{assumption}
                                              \label{assumption 10.3.1}
For any unit
$l\in\bR^{d}$ and $\eta>0$,  we have
$$
|\delta_{\eta,l}a^{\alpha}_{k}|\leq K_{0}(\sqrt{a^{\alpha}_{k}}+\eta).
$$

\end{assumption}

\begin{remark}
                                           \label{remark 10.3.1}
It is easy to see that
Assumption \ref{assumption 10.3.1} is satisfied
(with, perhaps, different $K_{0}$) if and only if
$\sigma^{\alpha}_{k}:=\sqrt{a^{\alpha}_{k}} $  
  satisfies the  Lipschitz
condition with constant $K_{0}$ with respect to $x$.

Indeed, the necessity follows after letting $\eta\downarrow
0$ and the sufficiency is a direct consequence of the formula
$$
 \delta_{\eta,l}a^{\alpha}_{k} = 2\sigma^{\alpha}_{k}
\delta_{\eta,l}\sigma^{\alpha}_{k}+\eta
(\delta_{\eta,l}\sigma^{\alpha}_{k})^{2} .
$$
Below we are also using the well-known fact that
a continuous function $v(x)$ is Lipschitz continuous 
with constant $K$ if and only if its generalized gradient
$v_{x}=D_{x}v $ satisfies $|v_{x}|\leq K$ (a.e.).
\end{remark}

\begin{definition}
                                            \label{definition 9.27.1}
Let $B$ be a finite subset of $\bR^{d}$
and
  $p(x,y)$ a real-valued function on $\bR^{d}\times\bR^{d}$.
For an $x_{0}\in\bR^{d}$
we say that the operator
\begin{equation}
                                                       \label{10.19.1}
S:u\to Su,\quad Su(x)=\sum_{y\in B\cup\{0\}}p(x,x+y)u(y)
\end{equation}
respects the maximum principle at $x_{0}$ relative to $B$
if, for any function $\phi(x)$ such that
$\phi(x_{0}+y)\geq\phi(x_{0})$
for all $y\in B$, we have $Su(x_{0})\geq0$.
\end{definition}

Obviously, the operators $\delta_{\eta, l}$ and $\Delta_{\eta, l}$
respect the maximum principle
at any point  relative to appropriate sets.

For $h>0$ set
\begin{equation}
                                                       \label{5.9.1}
 L^{\alpha}_{h}u=a^{\alpha}_{k}\Delta_{h,\ell_{k}} 
u +b^{\alpha}_{k}\delta_{h,\ell_{k}}u-c^{\alpha}u,
\end{equation}
where and throughout the paper the summation convention is enforced.
For each $t$ the operator $L^{\alpha}_{h}=L_{h}^{\alpha}(t,x)$ can be
considered as an operator on functions defined on $\bR^{d}$.
\begin{assumption}
                                              \label{assumption 9.27.3}
We have
$a^{\alpha}_{k}\geq h_{0}(b^{\alpha}_{k})^{-}$ (recall \eqref{5.13.1}).
\end{assumption}
\begin{remark}
                                       \label{remark 10.12.1}
It is easy to see that
Assumption \ref{assumption 9.27.3} implies that
for $h\in(0,h_{0}]$, $t\in\bR$, and $\alpha\in A$
the operator $L_{h}^{\alpha}(t,x)+c^{\alpha}(t,x)$
respects the maximum principle at any point $x_{0}$
relative to $ \Lambda_{0}$, where
\begin{equation}
                                                     \label{5.7.1}
\Lambda_{0}:=\{h \ell_{k}:k=\pm1,...,\pm d_{1}\}.
\end{equation}
In turn, provided that
all $\ell_{k}$ are different, 
the said property of $L_{h}^{\alpha}(t,x)+c^{\alpha}(t,x)$
implies what is required in 
Assumption \ref{assumption 9.27.3}.

\end{remark}

To satisfy Assumption \ref{assumption 9.27.3} it is sufficient
to require that $b^{\alpha}_{k}\geq0$, in which case 
what we use is just an 
  upwind
discretization of the ``transportation" term.

\begin{remark}
                                            \label{remark 10.19.1}
The operators $L^{\alpha}_{h}$ are natural
approximations of the operator
 \begin{equation}
                                                     \label{11.19.6}
L^{\alpha}u=a^{\alpha}_{k}\ell_{k}^{i}\ell_{k}^{j}u_{x^{i}x^{j}}
 +b ^{\alpha}_{k}\ell_{k}^{i}u_{x^{i}}-c u
\end{equation}
in the sense that $L_{h}^{\alpha}u\to L^{\alpha}u$ as
$h\downarrow0$ for all smooth $u$.

One may wonder how wide is the class of operators given
in the usual form
\begin{equation}
                                                 \label{5.7.2}
Lu=a^{ij}u_{x^{i}x^{j}}+b^{i}u_{x^{i}}
\end{equation}
which admit such a special approximation. 
We discuss this issue in Section \ref{section 5.7.1}.
\end{remark}

Next, we describe the free term in the equation, which are given 
by a real-valued function  
$$
f^{\alpha}
=f^{\alpha}(p,\psi,t,x)
$$ 
defined on $A\times\bR^{2d_{1}}\times\bR\times\bR\times\bR^{d}$.

\begin{assumption}
                                              \label{assumption 9.27.2}
The function $f^{\alpha}$ is bounded, $f^{\alpha}$ is continuous in
$\alpha$, continuous in $(p,\psi,x)$ and, for any $\alpha$
and $t$, its 
generalized gradients $D_{p}f^{\alpha}$, $D_{\psi}f^{\alpha}$,  and 
$D_{x}f^{\alpha}$ in $p,\psi$, and $x$, respectively, satisfy
$$
|D_{p^{k}}f^{\alpha}|\leq K_{0}\sqrt{a^{\alpha}_{k}},
\quad k=\pm1,...,\pm d_{1},\quad
|D_{\psi}f^{\alpha}|\leq K_{0},\quad
|D_{x}f^{\alpha}|\leq K_{3}
$$ 
for almost all $(p,\psi,x)\in\bR^{2d_{1}} \times\bR\times\bR^{d}$.
\end{assumption}

For fixed $h,\tau>0$ we consider the equation
\begin{equation}
                                                    \label{9.27.1}
\sup_{\alpha\in A}[r^{\alpha}
\delta^{T}_{\tau}u +L^{\alpha}_{h}u +
g^{\alpha} ]=0,
\end{equation}
where
$$
g^{\alpha}=g^{\alpha}(t,x)
=f^{\alpha}(\delta_{h,\ell_{k}}u(t,x),u(t,x),t,x).
$$
Observe that equation \eqref{9.27.1} takes the form
\eqref{10.26.1}. The presence of $r^{\alpha}$ in these equations
allows us to treat the normalized Bellman equations
(see \cite{Kr77}), which arise, for instance, in optimal
stopping problems or problems with singular control.

 Fix a vector $l\in\bR^{d}$
with $|l|\leq K_{0}$ and a number $\eta\in(0,h]$.
Set
$$
h_{i}=h\quad\text{for}\quad |i|=1,...,d_{1},\quad
h_{d+1}= h_{-(d+1)}=\eta,\quad\ell_{d+1}= -\ell_{-(d+1)}=  l,
$$
\begin{equation}
                                                    \label{5.11.1}
 \quad
\Lambda  =\{h_{1}\ell_{\pm1},...,h_{d+1}\ell_{\pm (d_{1}+1)}\}.
\end{equation}
We treat $\Lambda_{0}$ as a list rather than the set
with specified elements, even if $\ell_{1}=\ell_{2}$
we include in the list this vector  twice.

Observe that in \eqref{5.9.1} only $\ell_{k}\in\Lambda_{0}$
are involved. However, the method of ``shaking" the coefficients
requires estimates of difference derivatives in all directions
and not only along the mesh. This is the reason why
we introduce $\Lambda$. Set
$$
\Lambda_{n}=
\sum_{1}^{n}\Lambda=\{x:x=l_{1}+...+l_{n},\,\,l_{1},...,l_{n}\in\Lambda\},
\quad\Lambda_{\infty}=\bigcup_{n}\Lambda_{n}
$$
\begin{equation}
                                                      \label{5.9.4}
 \bar{\cM}_{T}=\{(n\tau)\wedge T:n=0,1,...\}\times\Lambda_{\infty},
\quad  \cM _{T}=\bar{\cM}_{T}\cap([0,T)\times\bR^{d}).
\end{equation}
 
Fix a  {\em finite\/} set $Q\subset \cM _{T}$,
assume that
\begin{equation}
                                                      \label{5.9.5}
Q|_{0}:=Q\cap(\{0\}\times\bR^{d})\ne\emptyset
\end{equation}
 and define
$$
\bar{Q}=Q\cup\{(t+\tau_{T}(t),x):(t,x)\in Q,t+\tau_{T}(t)=T\},
$$
\begin{equation}
                                                      \label{5.9.6}
 Q^{o}_{1}  =\{(t,x)\in Q: t<T, (t+\tau_{T}(t),x)\in \bar{Q},
(t,x +\Lambda)\subset Q  \},
\end{equation}
$$
 \partial_{1}  Q=\bar{Q}\setminus Q^{o}_{1}.
$$
Obviously, it may happen that $\bar{Q}=Q$. The subscript 1
is used above because later on we will need a ``fatter"
boundary $\partial_{2}Q$.

 Finally, define $T'$
as the least $n\tau$, $n=1,2,...$, such that $n\tau\geq T$,
recall that   $m\in\bR$ (see \eqref{5.13.1}) is a given fixed constant
 and 
introduce
$$
\xi (t)=e^{ m t},\quad t<T,\quad\xi (T)
=e^{ m T'},
\quad \xi_{(+)}=\xi\vee1,\quad\xi_{(-)}=\xi\wedge1 ,
$$
\begin{equation}
                                                      \label{5.9.7}
c_{m}=\frac{1-e^{-m\tau}}{\tau},\quad 
\lambda =
\inf_{\alpha,t,x}\big[
c^{\alpha}(t,x)+r^{\alpha}(t)c_{m} \big].
\end{equation}
Introducing a discontinuous function $\xi(t)$ may look
unnatural. However, what is important for us  is that 
\begin{equation}
                                                  \label{5.9.2}
\xi\delta_{\tau}^{T}u=e^{-m \tau}\delta_{\tau}^{T}(\xi u)
-c_{m}(\xi u)
\end{equation}
on $\cM_{T}$ for any $u=u(t,x)$.

 Everywhere below in this section $u$
is a given function on $\bar{\cM}_{T}$
satisfying \eqref{9.27.1} in $Q$. In our first result no  
control on the sizes of $r^{\alpha}$, $a^{\alpha}_{k}$,
$b^{\alpha}_{k}$,
and $c^{\alpha} $ is imposed.

\begin{theorem} 
                                            \label{theorem 9.16.1}
  Let    $h\in(0,h_{0}]$.
Then, under Assumptions \ref{assumption 9.27.1}
through \ref{assumption 9.27.2}
 there are constants $N=N(d_{1},K_{0})$,
  $N^{*}=N^{*}(d_{1},K_{0},K_{3})$ such that 
if $ \lambda \geq N$ then on $Q|_{0}$
\begin{equation}
                                               \label{9.16.2}
 |\delta_{\eta,l}u|\leq N^{*}e^{m^{+}(T+\tau)}\big[1+
\max_{\bar{Q}}|\xi_{(-)} u|
+\max_{k,\partial_{1} Q}(|\xi_{(-)}
\delta_{h,\ell_{k}}u|+|\xi_{(-)}\delta_{\eta,l}u|)
\big].
\end{equation}
 
\end{theorem}

We prove this theorem in Section \ref{section 10.3.1}.
\begin{remark}
                                              \label{remark 11.19.3}
This theorem is similar to Theorem 5.2 of \cite{Kr06}
and entails all the consequences derived from the latter
in \cite{Kr06} and \cite{DK3}. In particular, 
by using Theorem 5.6 of \cite{Kr06} and comparing
the equations for $u$ and $u(\tau+\cdot,\cdot)$
an estimate
of $\delta_{\tau} u$ can be obtained
if we require the data to have bounded derivatives in $t$.
\end{remark}   
 \begin{remark}
                                              \label{remark 5.9.1}
 Theorem \ref{theorem 9.16.1} has an immediate
application to elliptic equations. In that case
$u$ is independent of $t$, one can take $r^{\alpha}\equiv0$,
and use as large negative $m$ as one wishes without affecting
$\lambda$. Then it is seen that in \eqref{9.16.2}
the maximums over $\bar{Q}$ and
$\partial_{1}Q$ reduce to the maximums over $Q|_{0}$ and
$Q|_{0}\cap\partial_{1}Q$, respectively.

\end{remark}

Our second result is about Bellman
 equations with more general quasilinear operators.
This time \eqref{9.27.1} is assumed to be {\em uniformly
nondegenerate\/} in the space generated by $\ell_{k}$'s.
We will allow $f^{\alpha}(p,\psi,t,x)$ to grow   quadratically
with respect to $p$ and therefore no $b_{k}^{\alpha}$
are needed. The term $c^{\alpha}u$ also could be absorbed in
$f^{\alpha}$. However, we keep it, in order
  to state Theorem \ref{theorem 9.22.1}
in a simpler way.

\begin{assumption}
                                             \label{assumption 9.27.4}
 
(i) The functions $a^{\alpha}_{k}$ 
also depend on $\psi$: 
$$
a^{\alpha}_{k}=a^{\alpha}_{k}(\psi,t,x)
$$
and equation \eqref{10.26.1} holds in $Q $, where 
$F$ is defined by \eqref{11.9.1} with 
$a^{\alpha}_{k}(\psi,t,x)$ in place of $a^{\alpha}_{k}(t,x)$.
The functions $a^{\alpha}_{k}(\psi,t,x)$
are Lipschitz continuous
in $x$ with constant $K_{3}$, 
Lipschitz continuous
in $\psi$ with a constant $\omega\in(0,\infty)$,
$$
 a^{\alpha}_{k}\geq\delta,\quad |k|\leq
d_{1},\quad c^{\alpha}\geq-K_{3}.
$$

(ii)
The function $f^{\alpha}$ is continuous in $\alpha$,
continuous in $(p ,\psi , x )$, and
for   all values of the arguments,
satisfying $|\psi|\leq K_{1}$ and $|p|\geq K_{2}$, it holds that 
$$
 \quad |f^{\alpha}|  \leq \omega|p | ^{2 }
 +K_{3}.
$$

(iii) For each $\alpha$ and  $t$ 
the generalized gradients $D_{p}f^{\alpha}$,
$D_{\psi}f^{\alpha}$, and $D_{x}f^{\alpha}$ of $f^{\alpha}$
with respect to $p$, $\psi$, and $x$, respectively,
satisfy
$$
|D_{p}f^{\alpha}
|\leq \omega|p |   +K_{3} ,\quad 
| D_{\psi}f^{\alpha}|  \leq \omega  |p | ^{2 } 
 +K_{3}  ,
$$
\begin{equation}
                                                      \label{5.9.8}
|D_{x}f^{\alpha} |\leq \omega|p |  ^{3  }  
 +K_{3}
\end{equation}
(a.e.) on the set $\{|\psi|\leq K_{1}\}$.

 \end{assumption}
\begin{remark}
Clearly, Assumption \ref{assumption 9.27.4} (ii)
is satisfied with any $\omega>0$ 
and appropriate $K_{3}(\omega)$ if
$$
\sup_{\alpha,t,x,\psi}|f^{\alpha}(p,\psi,t,x)|=o(|p|^{2})
$$
as $|p|\to\infty$. This includes all functions affine in $p$
provided that the coefficients are bounded. Similar
situation occurs with 
Assumption~\ref{assumption 9.27.4}~(iii).
\end{remark}

\begin{assumption}
                                            \label{assumption 10.1.1}
For a constant $C\geq4$ depending only on $d_{1}$,
the exact value of which 
can be determined by examining the proof
of Theorem \ref{theorem 9.22.2}, we have
\begin{equation}   
                                           \label{10.1.3}
CK_{1}(1+K_{1})\omega\leq\delta.
\end{equation}
 
\end{assumption} 
 
\begin{theorem}
                                               \label{theorem 9.22.2}
Let 
$b^{\alpha}_{k} \equiv0$
and let Assumptions  \ref{assumption 9.27.1},
\ref{assumption 9.27.4}, 
and \ref{assumption 10.1.1}
 be satisfied. Assume that  $|u|\leq K_{1}$ in $\bar{Q}$ 
and $|\delta_{h,\ell_{k}}u| \leq K_{3}$
on $\partial_{1} Q$ if $|k|\leq d_{1}$. Then  
 in $\bar{Q}$  
$$
|\delta_{h,\ell_{k}}u|\leq N= N(d_{1},\delta, K_{1},K_{2},K_{3})
,\quad|k|\leq d_{1}.
$$
In particular, $N$ is independent of $T$.
\end{theorem}

The proof of this theorem is given in Section \ref{section 10.20.1}.

\begin{remark}
If $\omega$ is large, we need $K_{1}$ to be small in order to
satisfy \eqref{10.1.3}, that is, we need $u$ to be small.
By replacing $u$ with $u-\gamma$, where $\gamma$
is any constant, we see that, actually, we need
the oscillation of $u$ rather than $u$ itself
to be small if $\omega$ is
not. This restriction could be completely
avoided if we proved an  interior version 
of Theorem \ref{theorem 9.22.2} and a priori H\"older
continuity of $u$. It seems to the author that
this is possible, but requires much more work.
\end{remark}  

\begin{theorem} 
                                            \label{theorem 9.22.1}
Under the assumptions
of Theorem \ref{theorem 9.22.2} suppose
that  $ |\delta_{\eta,l}u|\leq K_{3}$
on $\partial_{1} Q$ and $a^{\alpha}_{k}$ are independent of
$\psi$. Then
there is a constant  $N=N(d_{1},\delta, K_{1},K_{2},K_{3})$,
   such that
if $ \lambda \geq N$,   then on $Q|_{0}$
$$
 |\delta_{\eta,l}u|\leq N e^{m^{+}(T+\tau)} .
$$
\end{theorem}

This is a simple corollary of Theorems \ref{theorem 9.22.2}
and \ref{theorem 9.16.1}
with $h_{0}=h$ in the latter. Indeed, once we know that
the values of $|\delta_{h,\ell_{k}}u|$ and $|u|$
are dominated by
a constant, the behavior of $f^{\alpha}(p,\psi,t,x)$
for large $|p|$ becomes irrelevant and we can even multiply it
by an appropriate cut-off function in such a way that
the new $f^{\alpha}$ would satisfy
Assumption \ref{assumption 9.27.2} and $u$ would still satisfy the new
equation.

Our next result is about second-difference estimates.

\begin{assumption}
                                              \label{assumption 10.2.2}
(i) The function $f^{\alpha}$ is independent of
$p$ and $\psi$.

(ii) For any $i,j=\pm1,...,\pm(d_{1}+1)$
and $\psi$ standing for any of the
functions   $b^{\alpha}_{k}$,
$c^{\alpha}$, and $f^{\alpha}$ 
we have
\begin{equation}
                                               \label{10.26.2}
|\delta_{h_{j},\ell_{j}}\delta_{h_{i},\ell_{i}}\psi|\leq K_{3},\quad  
|\delta_{h_{i},\ell_{i}}f^{\alpha}|\leq K_{3},\quad
|\delta_{h_{j},\ell_{j}}\delta_{h_{i},\ell_{i}}a^{\alpha}_{k}| 
\leq K_{0}+K_{3}\sqrt{a^{\alpha}_{k}}.
\end{equation}

\end{assumption}
Typical case when the third inequality in \eqref{10.26.2}
is satisfied  occurs if $a^{\alpha}_{k}=(\sigma^{\alpha}_{k})^{2}$,
where $\sigma^{\alpha}_{k}$ is bounded and twice continuously 
differentiable.

In contrast with the above results in which no control on the magnitudes
of $r^{\alpha}$, $a^{\alpha}_{k}$, $b^{\alpha}_{k}$, and
$c^{\alpha}$ is required, this time we need the following.

\begin{assumption}
                                              \label{assumption 10.3.2}
 We have
$$
\delta\leq\sup_{\alpha\in A}a_{k}^{\alpha} \leq K_{0},
\quad r^{\alpha}, |b^{\alpha}_{k}|,|c^{\alpha}|,|f^{\alpha}|\leq K_{3}.
$$
 
\end{assumption}

The following assumption is about a special structure
of the set of our basic vectors $\ell_{k},k=\pm1,...,\pm d_{1}$.
For $d_{1}=2$ and the standard grid
(generated by $\pm e_{1},\pm e_{2}$)
 it means that this set contains   all eight neighboring points of
the origin on the grid. 

\begin{assumption}
                                       \label{assumption 10.7.1}
There exists an integer $1\leq d_{0}<d_{1}$ such that
for  the list
\begin{equation}
                                                  \label{5.10.3}
  \cL :=\{h \ell_{\pm1},...,h\ell_{\pm d_{0}} \}
\end{equation}
and any $\ell_{k}$ with $d_{1}<|k|\leq d_{1}$ there exist
$l_{1},l_{2}\in\cL$ such that
$$
l_{1}\ne l_{2},\quad l_{1}\ne-l_{2},\quad \ell_{k}=l_{1}+l_{2}.
$$
 
\end{assumption}

One may think that Assumption \ref{assumption 10.7.1}
excludes the equations with only one spatial variable,
where it is natural to take $d_{1}=1$
and $\Lambda_{0}=\{\ell_{1},-\ell_{1}\}$.
However, we do not require
$\ell_{k}$ to be nonzero, and one can take $\Lambda_{0}$
to be $\{\ell_{1},-\ell_{1},\ell_{2},-\ell_{2},\ell_{3},
-\ell_{3}\}$ with $\ell_{2}=0$
and $\ell_{3}=\ell_{1}$. In that case 
Assumption \ref{assumption 10.7.1}
 is satisfied with $\cL=\{\ell_{1},
-\ell_{1},\ell_{2},-\ell_{2}\}$.
By the way the fact that now the origin is one of $\ell_{k}$
in no way contradicts Assumption \ref{assumption 10.3.2},
because in that case $\delta_{h,\ell_{k}}\phi=0$
and one can assign any value to $a^{\alpha}_{k}$
without changing the equation.

Define
$$
 Q^{o}_{2}  =\{(t,x)\in Q: t<T, (t+\tau_{T}(t),x)\in \bar{Q},
(t,x +\Lambda_{0}+\Lambda_{0})\subset Q  \},
$$
$$
 \partial_{2}  Q=Q\setminus Q^{o}_{2} .
$$
Here, naturally, $x +\Lambda_{0}+\Lambda_{0}=
\{x+y+z:y,z\in\Lambda_{0}\}$.
\begin{theorem}
                                              \label{theorem 10.4.1}
Suppose that Assumptions \ref{assumption 9.27.1},
\ref{assumption 10.3.1}, 
\ref{assumption 9.27.3},
\ref{assumption 10.2.2}-\ref{assumption 10.7.1}
are satisfied.
Then there exists a constant  $N=N(\delta, d_{1},K_{0})$
 such that
if $ \lambda \geq N$, then in $Q|_{0}$
for $i,j=\pm 1,...,\pm d_{1} $   we have
\begin{equation}
                                                 \label{10.8.1}
|\delta_{h,\ell_{j}}\delta_{h,\ell_{i}}u|
\leq 
N^{*}e^{m^{+}(T+\tau)}R , 
\end{equation}
where  $N^{*}=N^{*}(h_{0},\delta, d_{1},K_{0},K_{3})$,
$$
R = 1
+\max_{\bar{Q}}|\xi_{(-)}u| 
+\max_{Q}(\xi_{(-)}\delta_{\tau}^{T}u)^{-} 
$$
$$
+\max_{|i|,|j|\leq d_{1}}
\max_{\partial_{2}Q }
|\xi_{(-)}\delta_{h,\ell_{i} }\delta_{h,\ell_{j}}u|
+\max_{|i|\leq d_{1}}\max_{\bar{Q}}|\xi_{(-)}\delta_{h,\ell_{i}}u|).
$$
 
\end{theorem}
This theorem is proved in Section \ref{section 10.20.2}
 following a quite long
Section \ref{section 10.7.1} that contains the proof 
of Theorem \ref{theorem 10.4.1} under additional assumptions.

\begin{remark}
                                               \label{remark 11.19.4}
To get  ``closed" estimates of $\delta_{h,\ell_{j}}\delta_{h,\ell_{i}}u$
  we need to exclude $\delta_{h,\ell_{i}}u$ and
$\delta_{\tau}^{T}u$ from   $R$. This
can be done by using 
Theorem \ref{theorem 9.16.1} and
the idea from Remark \ref{remark 11.19.3}.
Another situation when 
$\delta_{\tau}^{T}u$ drops out
presents when $u$ is independent of $t$, so that, actually, we are
dealing with elliptic equations. We say more about this in the comments 
after Theorem \ref{theorem 11.2.1}
\end{remark}

In case of $a_{k}^{\alpha}$ independent of $x$
  Assumption \ref{assumption 10.7.1}
is not needed.

\begin{theorem}
                                              \label{theorem 11.2.1}
Suppose that Assumptions \ref{assumption 9.27.1},
\ref{assumption 10.3.1}, 
\ref{assumption 9.27.3},
\ref{assumption 10.2.2}, and \ref{assumption 10.3.2}
are satisfied. Also assume that $a^{\alpha}_{k}$
are independent of $x$,
$|\delta_{h_{i},\ell_{i}}b^{\alpha}_{k}|\leq K_{3}\sqrt{a^{\alpha}_{k}}$,
$|i|\leq d_{1}+1$, $|k|\leq d_{1}$, and $\lambda>0$.
Then   in $Q|_{0}$
for $k=\pm 1,...,\pm d_{1} $   we have
\begin{equation}
                                                    \label{11.4.4}
|\Delta_{h,\ell_{k}} u|
\leq 
N^{*}e^{m^{+}(T+\tau)}R_{0} , 
\quad (\Delta_{\eta,l}u)^{-}
\leq 
N^{*}e^{m^{+}(T+\tau)}R,
\end{equation}
where  $N^{*}=N^{*}(\lambda,h_{0},\delta, d_{1} ,K_{3})$,
$$
R_{0}= 1
+\max_{\bar{Q}}|\xi_{(-)}u| 
+\max_{Q }(\xi_{(-)}\delta_{\tau}^{T}u)^{-} 
$$
$$
+\max_{|i|\leq d_{1}}
\max_{\partial_{1} Q }|\xi_{(-)}\Delta_{h_{i},\ell_{i} } u|
+\max_{|i|\leq
d_{1}}\max_{\bar{Q}}|\xi_{(-)}\delta_{h_{i},\ell_{i}}u|) ,
$$
and $R$ is obtained from $R_{0}$ by taking $d_{1}+1$
in place of $d_{1}$.

\end{theorem}

This theorem  proved in Section \ref{section 11.4.1}
 is a direct generalization
of the corresponding result from \cite{DK1}: lower order coefficients
are allowed to depend on $(t,x)$
and we consider 
parabolic equations. In connection with the latter
observe that if $r^{\alpha}\equiv0$
(elliptic case), then one can 
let $m\to-\infty$ and see that 
in the definitions of $R_{0}$ and $R$ 
one can replace $\partial_{1} Q$ with
$Q|_{0}\cap\partial_{1} Q$.

\mysection{Some technical tools}

For any $\eta\in\bR^{d}$, $\nu\geq0$ set
 \begin{equation}                                  \label{5.11.2}
T_{\nu,\eta}\psi(x):=\psi(x+\nu\eta).
\end{equation}
\begin{lemma}
                                             \label{lemma 6.18.1}
For any  $\nu>0$,
 $l_{1},l_{2}\in\bR^{d}$, and  
functions $a(x),\psi(x)$ 
$$
\delta_{\nu ,l_{1}}(a\psi)=(\delta_{\nu ,l_{1}}a)\psi+
(T_{\nu ,l_{1}}a)\delta_{\nu ,l_{1}}\psi=a\delta_{\nu ,l_{1}}\psi+
\psi\delta_{\nu ,l_{1}}a
+\nu (\delta_{\nu ,l_{1}}a)\delta_{\nu ,l_{1}}\psi,
$$
$$
\delta_{\nu ,l_{2}}\delta_{\nu ,l_{1}}(a\psi)=
a\delta_{\nu ,l_{2}}\delta_{\nu ,l_{1}}\psi+
(\delta_{\nu ,l_{2}}a)\delta_{\nu ,l_{1}}\psi+
(\delta_{\nu ,l_{1}}a)\delta_{\nu ,l_{2}}\psi
$$
 \begin{equation}                                  \label{10.24.1}
+h[\delta_{\nu ,l_{1}}a+\delta_{\nu ,l_{2}}a]
\delta_{\nu ,l_{2}}\delta_{\nu ,l_{1}}\psi
+(\delta_{\nu ,l_{2}}\delta_{\nu,l_{1}}a)
T_{h,l_{1}+l_{2}}\psi,
\end{equation}
$$
\Delta_{\nu ,l_{1}}(a\psi)=a
\Delta_{\nu ,l_{1}}\psi+
\psi\Delta_{\nu ,l_{1}}a
+(\delta_{\nu ,l_{1}}a)\delta_{\nu ,l_{1}}\psi
+(\delta_{\nu ,-l_{1}}a)\delta_{\nu ,-l_{1}}\psi.
$$
In particular,
$$
\Delta_{\nu ,l_{1}}(\psi^{2})=2\psi\Delta_{\nu ,l_{1}}\psi
+(\delta_{\nu ,l_{1}}\psi)^{2}+(\delta_{\nu ,-l_{1}}\psi)^{2}.
$$
\end{lemma}

This lemma is proved by straightforward computations
(cf.~\cite{Kr06}). In the following lemma
we use Definition \ref{definition 9.27.1}.

\begin{lemma}
                                               \label{lemma 6.18.2}
 
If an operator
$$
S\psi(x)=\sum_{y\in\Lambda\cup\{0\}}p(x,x+y) \psi(x+y) 
$$
respects the maximum principle at a point $x_{0}\in\bR^{d}$
relative to $ \Lambda $ and $\psi$ is a function such that
$\psi(x_{0})\leq0$, then  
$
-S\psi\leq S(\psi^{-})$ at $x_{0}$. In particular,
$\phi^{-}S\phi^{-}\geq-\phi^{-}S\phi$ at $x_{0}$
for any function $\phi$.

\end{lemma}

This follows from the definition
and the fact that $\psi+
\psi^{-}\geq0$ on $x_{0}+\Lambda$ and $\psi+\psi^{-}=0$ at
$x_{0}$.

 The following lemma from \cite{Kr06} is used in the proof of
Theorem~\ref{theorem 10.4.1}.
\begin{lemma}
                                               \label{lemma 6.18.02}

Let $\psi$ be a function on $\bR^{d}$, $\nu>0$.
Then  
\begin{equation}
                                          \label{6.18.6}
|\Delta_{\nu,\eta}\psi|\leq 
|\delta_{\nu,-\eta}((\delta_{\nu,\eta}\psi)^{-})|+
 |\delta_{\nu,\eta}((\delta_{\nu,-\eta}\psi )^{-})| .
\end{equation}

\end{lemma}

\mysection{Proof of Theorem \protect\ref{theorem 9.16.1}}
                                              \label{section 10.3.1}

We start with some preparations.
From now on index $k$ will run through $\{\pm1,...,\pm d_{1}\}$
and $i,j$ through $\{\pm1,...,\pm (d_{1}+1)\}$.
By $N$ and $N^{*}$ in this section we denote
generic constants depending on the data
as in the statement of the theorem.
We use the notation \eqref{5.9.4} through \eqref{5.9.7}
and introduce few new objects.
We need two  constants $\varepsilon$ and $\mu$ defined by
$$
\varepsilon^{-1}-2\varepsilon d_{1}=1,\quad4\mu =
(d_{1}+1)^{-1}\wedge\varepsilon.
$$

Introduce
$$
\Gamma =\{\gamma=(\gamma_{i}:i=\pm1,...,\pm (d_{1}+1)):
\gamma_{i}\in[\varepsilon,\varepsilon^{-1}]\},
$$
$$
\delta_{i}=\delta_{h_{i},\ell_{i}},\quad
P_{\gamma}\phi=\gamma_{i}\delta_{i}\phi,
\quad v=\xi u,\quad\Delta_{k}=\Delta_{h_{k},\ell_{k}}.
$$

By using \eqref{5.9.2} we see that in $Q$
$$
F\big(
e^{-m \tau}\delta_{\tau}^{T}v
-c_{m}v, \xi\Delta_{k}u,
\xi\delta_{k}u,\xi u\big)=0.
$$

Also introduce
$$
v_{\gamma}= P_{\gamma}v,\quad
v_{i}=\delta_{i} v,
$$
$$
P_{\gamma\mu}\phi =v_{\gamma}^{-}
P_{\gamma}\phi 
-\mu v_{i}
 \delta_{i}  \phi ,
$$
$$
 W=
\sum_{i} v_{i} ^{2},
\quad
V_{\gamma\mu}=[v_{\gamma}^{-}]^{2}
+\mu W.
$$
Observe that 
$$
P_{\gamma\mu}v=-V_{\gamma\mu}.
$$

Finally, 
let $(\gamma_{0},t_{0},x_{0})\in\Gamma \times\bar{Q}$ be a
point at which $V_{\gamma\mu}$ attains its maximum value over
$\Gamma \times\bar{Q}$.  
 \begin{theorem}
                                                \label{theorem 10.11.1}
The assertions of Theorem \ref{theorem 9.16.1}
hold  true   if in addition to its assumptions
$(t_{0},x_{0})\in Q^{o}_{1}$ and
 \begin{equation}
                                                          \label{9.8.04}
v_{ \gamma_{0}}^{-}(t_{0},x_{0})\geq(1/2)\max_{\bar{Q},i}|v_{i}|.
\end{equation}
\end{theorem}

To prove this theorem we need an auxiliary result.

\begin{lemma}
                                              \label{lemma 10.11.3}
Assume \eqref{9.8.04}. Then
the operator $P_{ \gamma_{0}\mu }$
respects the maximum principle at $(t_{0},x_{0})$, that is,
for any   function $\phi$ such that
 $\phi(x_{0})\geq\phi(x_{0}+ \eta)$ for all
$\eta\in\Lambda$, we have $P_{\gamma_{0}\mu}\phi(t_{0},x_{0})
\leq0$.

\end{lemma}

Proof. Since $P_{\gamma\mu}1=0$ we may assume that
$\phi(x_{0})=0$. Then at $(t_{0},x_{0})$
$$
P_{\gamma_{0}\mu}\phi =(v_{\gamma_{0}}^{-}\gamma_{0i}
-\mu v_{i})\phi(x_{0}+h_{i}\ell_{i}),
$$
which is negative since $\phi(x_{0}+h_{i}\ell_{i})\leq0$
and $v_{\gamma_{0}}^{-}\gamma_{0i}
-\mu v_{i}\geq v_{\gamma_{0}}^{-}(\varepsilon-2\mu)\geq0$
($\mu\leq\varepsilon/2$).
The lemma is proved.

We also need the following construction.
Notice that, if $(t_{0},x_{0})\in Q$,
 there is a sequence $\alpha_{n}\in A$
such that at $(t_{0},x_{0})$
$$
\lim_{n\to\infty}[
r^{\alpha_{n}} \delta^{T}_{\tau}u +
a^{\alpha_{n}}_{k} 
\Delta_{k}u 
+b^{\alpha_{n}}_{k} 
\delta_{k}u 
-c^{\alpha_{n}} u +
g^{\alpha_{n}} ]
$$
$$
= \sup_{\alpha\in A}[r^{\alpha} 
\delta^{T}_{\tau}u +
a^{\alpha}_{k}
\Delta_{k}u +b^{\alpha}_{k}
\delta_{k}u 
-c^{\alpha}u+g^{\alpha}] =0.
$$
Since the numbers of possible values of $t$ for points in $Q$ is finite,
and the functions 
$a^{\alpha}_{k}
(t,x)$, $b^{\alpha}_{k}(t,x)$, $c^{\alpha}(t,x)$, $f^{\alpha}(p,\psi,t,x)$
are uniformly continuous functions of $(p,\psi, x)$,
there is a subsequence $\{n'\}\subset\{1,2,...\}$
and  functions $\bar{r}$, 
$\bar{a}_{k}
(t,x)$, $\bar{b}_{k}(t,x)$, $\bar{c}(t,x)$, $\bar{f}(p,\psi,t,x)$
such that they satisfy our
assumptions  changed
in an obvious way and
$$
(\bar{r},\bar{a}_{k}
(t,x),\bar{b}_{k}(t,x),\bar{c}(t,x),\bar{f}(p,\psi,t,x))
$$
$$
=\lim_{n'\to\infty}
(r^{\alpha_{n'}}(t),a^{\alpha_{n'}}_{k}(t,x),b^{\alpha_{n'}}_{k}(t,x),
c^{\alpha_{n'}} (t,x),f^{\alpha_{n'}} (p,\psi,t,x))
$$
on $Q$ for all $p,\psi$.  

Obviously,  
for
$$
\bar{g}(t,x):=\bar{f}(\delta_{k}u(t,x),u(t,x),t,x)
$$
at $(t_{0},x_{0})$ we have
$$
\bar{r}\delta^{T}_{\tau}u+\bar{a}_{k}\Delta_{k}
u+\bar{b}_{k}
\delta_{k}u-\bar{c}u+\bar{g}=0,
$$
\begin{equation}
                                               \label{9.7.01}
e^{-m \tau}\bar{r}\delta^{T}_{\tau}v
+ \bar{a}_{k}\Delta_{k}
v+\bar{b}_{k}
\delta_{k}v-(\bar{c}+\bar{r}c_{m})v+\xi\bar{g}=0
\end{equation}
and, if $(t_{0},x_{0})\in Q^{o}_{1}$, then
for any $i$ ($ =\pm1,...,\pm (d_{1}+1)$, the shift operator
$T$ is introduced in \eqref{5.11.2}) 
\begin{equation}
                                               \label{9.7.02}
T_{h_{i },\ell_{i }}[e^{-m \tau}\bar{r}\delta^{T}_{\tau}v
 +
\bar{a}_{k}\Delta_{k}
v+\bar{b}_{k}
\delta_{k}v-(\bar{c}+\bar{r}c_{m})v+\xi\bar{g}]\leq0,
\end{equation}
where and below 
for simplicity of notation
we drop $(t_{0},x_{0})$ in the arguments of functions
  we are dealing with. 

{\bf Proof of Theorem \ref{theorem 10.11.1}}.  
Set 
$$
\bar{L}^{0}_{h} =\bar{a}_{k}\Delta_{k}
+\bar{b}_{k}\delta_{k},\quad 
\bar{L} _{h}=\bar{L}^{0}_{h}-(\bar{c}+\bar{r}c_{m}).
$$
By Assumption \ref{assumption 9.27.3},
 Lemma   \ref{lemma 6.18.1}, and
 Lemma  \ref{lemma 6.18.2} (with $v_{\gamma }$ in place of~$\phi$) 
$$
0\geq\bar{L}^{0}_{h}V_{\gamma_{0}\mu} 
=2v_{\gamma_{0}}^{-}\bar{L}_{h}^{0}v^{-}_{\gamma_{0}} 
+2\mu v_{i}\bar{L}^{0}_{h}v_{i}+ I_{1}+ \mu I_{2} 
$$
$$
\geq -2v_{\gamma_{0}}^{-}\bar{L}_{h}^{0}v_{\gamma_{0}} 
+2\mu v_{i}\bar{L}^{0}_{h}v_{i}+ I_{1}+ \mu I_{2},
$$
where
$$
I_{1}:=2\bar{a}_{k} (\delta_{k}v_{\gamma_{0}}^{-}) ^{2}
+h\bar{b}_{k}(\delta_{k}v_{\gamma_{0}}^{-})^{2},
$$
$$
I_{2}:=
2\bar{a}_{k} \sum_{i}(\delta_{k}v_{i} ) ^{2}
+h\bar{b}_{k}\sum_{i}(\delta_{k}v_{i} )^{2}.
$$
Since $2\bar{a}_{k}+h\bar{b}_{k}\geq\bar{a}_{k}$ we conclude
\begin{equation}
                                                        \label{9.9.5}
v_{\gamma_{0}}^{-}\bar{L}_{h}^{0}v_{\gamma_{0}} 
-\mu v_{i}\bar{L}^{0}_{h}v_{i}
\geq  (1/2)\mu\bar{a}_{k} \sum_{i}(\delta_{k}v_{i} ) ^{2}.
\end{equation}

On the other hand, by \eqref{9.7.01} and \eqref{9.7.02}
and Lemma \ref{lemma 10.11.3}
 at $(t_{0},x_{0})$
\begin{equation}
                                                        \label{9.13.1}
  P_{\gamma_{0}\mu}
[e^{-m \tau}\bar{r}\delta^{T}_{\tau}v
+ 
 \bar{L}_{h}v+\xi\bar{g}]\leq0.
\end{equation}
Owing to \eqref{9.9.5} we obtain
$$
P_{\gamma_{0}\mu}\bar{L}_{h}v
=v_{\gamma_{0}}^{-}\bar{L}_{h}^{0}v_{\gamma_{0}}
-\mu v_{i}\bar{L}_{h}^{0}v_{i}+
(\bar{c}+\bar{r}c_{m})V_{\gamma_{0}\mu}+I_{3}+I_{4}+I_{5}
$$
$$
\geq \lambda V_{\gamma_{0}\mu}
+   
(1/2)\mu\bar{a}_{k}
\sum_{i}(\delta_{k}v_{i} ) ^{2} +I_{3}+I_{4}+I_{5},
$$
where
$$
I_{3}:=v_{\gamma_{0}}^{-}
(P_{\gamma_{0}}\bar{a}_{k})\Delta_{k}
v+h_{i}v_{\gamma_{0}}^{-}
 \gamma_{0i}(\delta_{i}\bar{a}_{k})\Delta_{k}v_{i}
$$
$$
-\mu v_{i}(\delta_{i}\bar{a}_{k})\Delta_{k}v
-\mu
h_{i}v_{i}(\delta_{i}\bar{a}_{k})\Delta_{k}v_{i},
$$
$$
I_{4}:=[v_{\gamma_{0}}^{-}\gamma_{0i}
  -\mu
v_{i} ](\delta_{i}\bar{b}_{k} )T_{h_{i},\ell_{i}} v_{k},
$$
$$
I_{5}: =[v_{\gamma_{0}}^{-}
 \gamma_{0i}    
-\mu v_{i}  ]
(\delta_{i}\bar{c} ) T_{h_{i},\ell_{i}}v.
$$

Upon observing that 
$$
\Delta_{k}
=-\delta_{k}\delta_{-k},\quad 
 h\Delta_{k}=\delta_{k} +\delta_{-k} 
,\quad h_{i}\leq h
$$
and by assumption \eqref{9.8.04}
$$
h|\delta_{k}v_{-k}|\leq2
\max_{\bar{Q},i}|v_{i}|\leq4v_{\gamma_{0}}^{-},
\quad
|h^{2}\Delta_{k}v_{i}|\leq4\max_{\bar{Q}}|v_{i}|
\leq 8v_{\gamma_{0} }^{-},
$$
 we find
$$
|I_{3}|\leq
Nv_{\gamma_{0} }^{-}
(\sqrt{\bar{a}_{k}}+h)|\delta_{k}v_{-k}|
+Nv_{\gamma_{0}}^{-}
\sum_{i}(\sqrt{\bar{a}_{k}}+h)|h\Delta_{k}v_{i}|
$$
$$
\leq
Nv_{\gamma_{0} }^{-} \sum_{i}\sqrt{\bar{a}_{k}}
|\delta_{k}v_{i}| +N(v_{\gamma_{0} }^{-})^{2}\leq 
N(v_{\gamma_{0} }^{-})^{2}+
(1/4)\mu\bar{a}_{k} \sum_{i}(\delta_{k}v_{i} ) ^{2}.
$$
 
This   yields
$$
P_{\gamma_{0}\mu}\bar{L}_{h}v
\geq 
  \lambda V_{\gamma_{0}\mu} -N (v_{\gamma_{0}
}^{-})^{2}+(1/4)\mu\bar{a}_{k}
\sum_{i}(\delta_{k}v_{i} ) ^{2} +I_{4}+I_{5}.
$$

Next, obviously,
$$
|I_{4}|\leq N(v_{\gamma_{0} }^{-})^{2},\quad
|I_{5}|\leq N^{*} v_{\gamma_{0} }^{-} \max_{\bar{Q}}|v|.
$$
Therefore, and since $(v_{\gamma_{0} }^{-})^{2}\leq V_{\gamma_{0}\mu}$,
\begin{equation}
                                                \label{9.13.2}
P_{\gamma_{0}\mu}\bar{L}_{h}v
\geq 
 (\lambda -N)V_{\gamma_{0}\mu}
-N^{*}\max_{Q}|v|^{2}
+(1/4)\mu\bar{a}_{k} \sum_{i}(\delta_{k}v_{i} ) ^{2}.
\end{equation}

Now we deal with other terms in \eqref{9.13.1}.
Note that, since $(t_{0}+\tau_{T}(t_{0}),x_{0})\in\bar{Q}$,  
$$
0\geq \delta^{T}_{\tau}V_{\gamma_{0}\mu}
=2v^{-}_{\gamma_{0}}\delta^{T}_{\tau}
v^{-}_{\gamma_{0}}+2\mu v_{i}\delta^{T}_{\tau}v_{i}
+\tau (\delta^{T}_{\tau}
v^{-}_{\gamma_{0}})^{2}+\mu\tau \sum_{i}(\delta^{T}_{\tau}v_{i})^{2},
$$
so that (cf.~Lemma \ref{lemma 6.18.2})
$$
0\geq-v^{-}_{\gamma_{0}}\delta^{T}_{\tau}
v _{\gamma_{0}}+\mu v_{i}\delta^{T}_{\tau}v_{i}=-P_{\gamma_{0}\mu}
\delta^{T}_{\tau}v.
$$
By recalling that $r^{\alpha}\geq0$,
we find from \eqref{9.13.2} and \eqref{9.13.1} that
\begin{equation}
                                                      \label{9.27.2}
 N^{*}\max_{\bar{Q}}|v|^{2}\geq (  \lambda  -N)V_{\gamma_{0}\mu}
+\xi  P_{\gamma_{0}\mu}\bar{g}
+ (1/4)\mu\bar{a}_{k} \sum_{i}(\delta_{k}v_{i} ) ^{2} .
\end{equation}

By Assumption \ref{assumption 9.27.2} 
$$
|\xi P_{\gamma_{0}\mu}\bar{g}|\leq N V_{\gamma_{0}\mu}^{1/2}
(\sum_{i}\sqrt{\bar{a}_{k}}|\delta_{k}v_{i}|
+V_{\gamma_{0}\mu}^{1/2}+N^{*}\xi ).
$$
It follows that the sum of the last two terms in \eqref{9.27.2}
is greater than
$$
-NV_{\gamma_{0}\mu}- N^{*}\xi V_{\gamma_{0}\mu}^{1/2}
\geq
-N V_{\gamma_{0}\mu}-N^{*}\xi^{2}
$$ 
Hence
$$
N^{*}(\max_{\bar{Q}}|v|^{2}+\xi^{2}) 
\geq (  \lambda  -N_{1})V_{\gamma_{0}\mu},
$$
which for $\lambda-N_{1}\geq1$ shows  
that on $Q|_{0}$
$$
 |\delta_{\eta,l}u|^{2}=
 |\delta_{\eta,l}v|^{2}\leq\mu^{-1 }V_{\gamma_{0}\mu}
$$
$$
\leq N^{*}(\max_{\bar{Q}}|v|^{2}+\xi^{2})
\leq N^{*}\xi_{(+)}^{2}(T)(\max_{\bar{Q}}|\xi_{(-)}u|^{2}+1).
$$
This implies \eqref{9.16.2} and the theorem is proved.

 In light of this theorem
to prove Theorem \ref{theorem 9.16.1} we only need
to show that in case  the assumption of 
Theorem \ref{theorem 10.11.1}
are not satisfied one can obtain the assertion
of Theorem \ref{theorem 9.16.1} differently.
To do that we need two lemmas.

\begin{lemma}
                                               \label{lemma 9.8.1}
Take a  function $\phi$ on $\bR^{d}$ and assume that
$$
\max_{i}|\delta_{i}\phi(x)|
\leq\max_{i}|\delta_{i}\phi(0)|
$$
for $x\in  \Lambda$. Then
\begin{equation}
                                               \label{9.8.1}
\max_{i}|\delta_{i}\phi(0)|\leq  
\max\big((P_{\gamma}\phi(x))^{-}:
\gamma\in \Gamma ,x\in\{0\}\cup  \Lambda\big).
\end{equation}

\end{lemma}

Proof. 
 Take
a  $j$ such that
$$
\max_{i}|\delta_{i}\phi(0)|=
|\delta_{j}\phi(0)|
$$
and first assume that
\begin{equation}
                                               \label{9.8.2}
|\delta_{j}\phi(0)|
=-\delta_{j}\phi(0).
\end{equation}
Then 
take $\gamma_{j}=\varepsilon^{-1}$ and $\gamma_{i}=\varepsilon$
for $i\ne j$.
Since $-\delta_{j}\phi(0)\geq
\delta_{i}\phi(0)$
and $-\delta_{j}\phi(0)\geq0$, we have
$$
P_{\gamma}\phi(0)=\varepsilon^{-1}
\delta_{j}\phi(0)
+\varepsilon\sum_{i\ne j}\delta_{i}\phi(0)
$$
$$
\leq
(\varepsilon^{-1}-2\varepsilon
d_{1})\delta_{j}\phi(0)
=-(\varepsilon^{-1}-2\varepsilon d_{1})
\max_{i}|\delta_{i}\phi(0)|
$$
and \eqref{9.8.1} follows since $\varepsilon^{-1}-2\varepsilon d_{1}=1$.

If \eqref{9.8.2} is not satisfied, then
$$
\max_{i}|\delta_{i}\phi(h_{j}\ell_{j})|\leq
|\delta_{j}\phi(0)|
=\delta_{j}\phi(0)=-
\delta_{-j}\phi(h_{j}\ell_{j}),
$$
which combined with an obvious inequality 
between the extreme terms
yields
$$
\max_{i}|\delta_{i}\phi(h_{j}\ell_{j})|=
-
\delta_{-j}\phi(h_{j}\ell_{j})
=|
\delta_{-j}\phi(h_{j}\ell_{j})|
$$
$$
=|\delta_{j}\phi(0)|=
\max_{i}|\delta_{i}\phi(0)|.
$$
By the above argument applied to the point $h_{j}\ell_{j}$
in place of $0$
$$
\min_{\Gamma }P_{\gamma}\phi(h_{j}\ell_{j})\leq
(\varepsilon^{-1}-2\varepsilon d_{1})\delta_{-j}
\phi(h_{j}\ell_{j})
$$
$$
=-(\varepsilon^{-1}-2\varepsilon d_{1})
\max_{i}|\delta_{i}\phi(0)|=
-\max_{i}|\delta_{i}\phi(0)|.
$$
  The lemma is proved.

\begin{lemma}
                                           \label{lemma 9.8.3}
 
Condition \eqref{9.8.04} is satisfied
if $(t_{0},x_{0})\in Q^{o}_{1}$ and
\begin{equation}
                                                         \label{10.12.2}
\max_{\bar{Q},i}v_{i}^{2}=\max_{Q^{o}_{1},i}v_{i}^{2}.
\end{equation}

\end{lemma}

Indeed if \eqref{9.8.04} does not  hold, then
$$
\max_{\Gamma\times\bar{Q}}[v_{\gamma}^{-}]^{2}
\leq V_{\gamma_{0}\mu}(t_{0},x_{0})
$$
$$
<(2(d_{1}+1)\mu+1/4)\max_{\bar{Q},i}v_{i}^{2}
=(2(d_{1}+1)\mu+1/4)\max_{Q^{o}_{1},i}v_{i}^{2},
$$
which contradicts Lemma \ref{lemma 9.8.1}
since $ \mu \leq1/(4(d_{1}+1))$.
 
{\bf Proof of Theorem \ref{theorem 9.16.1}}.
If $(t_{0},x_{0})\in
\partial_{1} Q$, then in $Q|_{0}$
$$
 |\delta_{\eta,l}u|= |\delta_{\eta,l}v|
\leq\mu^{-1/2}V_{\gamma_{0}\mu}^{1/2}\leq
\mu^{-1/2}V_{\gamma_{0}\mu}^{1/2}(t_{0},x_{0}),
$$
where the last term is obviously less than
the right-hand side of \eqref{9.16.2}. 
Furthermore, if $(t_{0},x_{0})\in Q^{o}_{1}$ but
\eqref{9.8.04} is violated, then by Lemma \ref{lemma 9.8.3}
$$
\max_{\bar{Q} ,i}v_{i}^{2}=
\max_{\partial_{1} Q,i}v_{i}^{2}
$$
and one can use the above argument. 
Finally,
in the remaining case
Theorem \ref{theorem 10.11.1} is applicable.
The theorem is proved.

\begin{remark}
As in \cite{DK1} one can relax conditions on
$\lambda$ by adding in its expression any large constant times
$\min_{k}a^{\alpha}_{k}$. This can be shown
by considering maximum points of $V_{\gamma\mu}+\nu v^{2}$,
where $\nu$ is a large constant. This would also allow us
to relax the condition on $D_{\psi}f$ to
$|D_{\psi}f|\leq K_{0}+K_{3}\min_{k}a^{\alpha}_{k}$.
\end{remark}

\mysection{Proof of Theorem \protect\ref{theorem 9.22.2}}
                                             \label{section 10.20.1}

Our goal is to show how to choose an appropriate
$C=C(d_{1})$ in \eqref{10.1.3}. Below in this section
by $N$ we denote   generic constants depending only
on $d_{1},\delta,  K_{1},K_{2}, K_{3}$ but not $\omega$. 

First of all, observe that $l$ does not enter  either
equation \eqref{9.27.1}  or the statement of the theorem.
It is involved, however, in the definition of 
$\partial_{1} Q$ making it ``fatter". Because of that if we
additionally assume that $l=0$,  the result will be stronger.
Therefore,   we assume that $l=0$.

Set $m =0$ 
 and introduce $\varepsilon$,  
$\Gamma $, $v_{\gamma}$, $v_{i}$,
$V_{\gamma\mu}$, $W$,
$P_{\gamma}$, $P_{\gamma\mu}$ as 
in the beginning of Section \ref{section 10.3.1}.
However, since $v=u$,
we also write $u_{\gamma}$,  $u_{i} $, 
and $U_{\gamma\mu}$ instead
$ v_{\gamma}$, $ v_{i}$, and $V_{\gamma\mu}$, respectively.
Since $l=0$, $u_{\pm (d_{1}+1)}=v_{\pm(d_{1}+1)}=0$
and now there is no need to allow $i$ to take the values
$\pm(d_{1}+1)$. Therefore, we restrict it to
the range $\pm1,...,\pm d_{1}$.
As everywhere in the article index $k$ runs through
$\pm1,...,\pm d_{1}$.
 This time we take
$$
8\mu=d_{1}^{-1}\wedge\varepsilon.
$$

Introduce $\kappa \geq0$ as a solution of
\begin{equation}
                                           \label{9.21.1}
\max_{\Gamma \times\bar{Q}}[U_{\gamma\mu}
+\kappa  u^{2}] =4\kappa  K_{1}^{2} .
\end{equation}
Observe that if $\kappa=0$ is a solution of \eqref{9.21.1},
then $u_{k}\equiv0$ and the assertion of the theorem
is trivial. Therefore, without losing generality we may assume
that for $\kappa=0$ the left-hand side of \eqref{9.21.1}
is strictly greater than its right-hand side. Furthermore,
as a function of $\kappa$ the left-hand side is convex increasing
with Lipschitz constant not greater than $K_{1}^{2}$.
It follows that \eqref{9.21.1} has a unique solution $\kappa>0$.
After that we define $(\gamma_{0},t_{0},x_{0})$ as a point
in $\Gamma \times \bar{Q}$ at which
the maximum in \eqref{9.21.1} is attained.
 For simplicity of notation we drop the arguments
$(t_{0},x_{0})$ in what follows. We also use
the abbreviated notation $\delta_{i},\Delta_{k}$
introduced in Section \ref{section 10.3.1}.

First we show that $
[u_{\gamma_{0}}^{-} ]^{2}$ is the main term in $U_{\gamma_{0}\mu}
+\kappa u^{2}$.
\begin{lemma}
                                           \label{lemma 9.21.1}
Assume that \eqref{10.12.2} holds. Then
\begin{equation}
                                           \label{9.21.2}
[u_{\gamma_{0}}^{-} ]^{2}\geq 
2\kappa  K_{1}^{2} ,
\end{equation}
\begin{equation}
                                           \label{9.21.3}
2d _{1}\varepsilon^{-1}\max_{i}|u_{i} | \geq
 u_{\gamma_{0}}^{-} \geq (1/2) \max_{\bar{Q},i} |u_{i}| .
\end{equation}
In particular $u_{\gamma_{0}}  <0$
and   the operator
$P_{\gamma_{0}\mu}$
respects the maximum principle at $(t_{0},x_{0})$.
\end{lemma}

Proof. First, notice that, owing to \eqref{10.12.2},
 Lemma \ref{lemma 9.8.1}, and the definitions
of $\kappa$ and $(t_{0},x_{0})$
we have  
$$
4\kappa K_{1}^{2} =[u_{\gamma_{0}}^{-}]^{2}+\mu\sum_{i}u_{i}^{2}
+\kappa u^{2}
\leq
[u_{\gamma_{0}}^{-}]^{2}+2 d_{1} \mu
\max_{\bar{Q},i} u_{i} ^{2}
+\kappa K_{1}^{2} 
$$
$$
\leq
[u_{\gamma_{0}}^{-}]^{2}+2 d_{1} \mu
\max_{\Gamma\times\bar{Q}}[u_{\gamma }^{-}]^{2}
+\kappa K_{1}^{2} 
\leq [u_{\gamma_{0}}^{-}]^{2}+2 d_{1} \mu 4
\kappa K_{1}^{2} +\kappa K_{1}^{2}.
$$
This implies \eqref{9.21.2} since $8 d_{1} \mu\leq1$.

The first inequality in \eqref{9.21.3} is obvious. If the second one
is wrong, then  
$$
\max_{\Gamma \times \bar{Q}}[u_{\gamma}^{-}]^{2}
\leq U_{\gamma_{0}\mu} +\kappa 
u^{2} 
< [2 d_{1} \mu+1/4]\max_{\bar{Q},i}u_{i}^{2}+\kappa K_{1}^{2} 
$$
$$
\leq [2 d_{1} \mu+1/4]\max_{\Gamma \times
\bar{Q}}[u_{\gamma}^{-}]^{2}+\kappa K_{1}^{2}
\leq (1/2)\max_{\Gamma \times
\bar{Q}}[u_{\gamma}^{-}]^{2}+\kappa K_{1}^{2} 
$$
 contrary to  \eqref{9.21.2}.
  The last assertion of the lemma
follows from Lemma \ref{lemma 10.11.3}. The lemma is proved.

Next, if $(t_{0},x_{0})\in\partial_{1} Q$, then
$$
4\kappa K_{1}^{2} =
 U_{\gamma_{0}\mu} 
+\kappa  u^{2}  \leq N'+\kappa
K_{1}^{2} ,\quad
3\kappa K_{1}^{2} \leq N',
$$
and $u_{k}^{2}(t,x)\leq (4/3)\mu^{-1}N'$ in $\bar{Q}$ according to
\eqref{9.21.1}. In this case the assertion of the theorem is true.
Similarly, as in the  proof of Theorem~\ref{theorem 9.16.1}
we get the result if \eqref{10.12.2} is violated.
This justifies the first two assumptions in
the following set   which we impose:
$$
(t_{0},x_{0})\in Q^{o}_{1},\quad
\max_{\bar{Q},k}u^{2}_{k}=\max_{Q^{o}_{1},k}u^{2}_{k},
$$
\begin{equation}
                                                   \label{9.23_1}
4d_{1}\varepsilon^{-1} W^{1/2}\geq 
\max_{Q^{o}_{1},k}|u _{k}|,\quad
W \geq K_{2}^{2}+ 8K_{1} K_{3}\delta^{-1}+1.
\end{equation}
The third relation in \eqref{9.23_1} follows
from Lemma \ref{lemma 9.21.1} and
the last assumption in \eqref{9.23_1} restricts us
 to the only
nontrivial case in light of~\eqref{9.21.3}.

Next again owing to the fact that the number of points
in $Q$ is finite we can find some functions $\bar{r}(t)$,
$\bar{a}_{k}( t,x)$, $\bar{c}( t,x)$, and $\bar{f}(p,\psi,r,t,x)$
satisfying Assumptions \ref{assumption 9.27.1} and 
\ref{assumption 9.27.4} and such that for
$$
 \bar{g}(t,x):=\bar{f}(u_{k}(t,x),u(t,x),t,x)-
 \bar{c} (t,x)u(t,x),
\quad 
\bar{L}_{h}:=\bar{a}_{k}\Delta_{k}
$$
we have (recall that the arguments $(t_{0},x_{0})$
are dropped)
\begin{equation}
                                           \label{9.21.4}
\bar{r}\delta^{T}_{\tau}u+\bar{L}_{h}
u+\bar{g}=0,\quad 
T_{h ,\ell_{i }}[\bar{r}\delta^{T}_{\tau}u+
\bar{L}_{h}u+\bar{g} ]\leq0
\end{equation}
for any $i$ ($ =\pm1,...,\pm  d_{1}  $). 

Now, as before  
$$
0\geq \bar{L}_{h}[U_{\gamma_{0}\mu}
+\kappa u^{2} ]\geq
-2u_{\gamma_{0}}^{-}\bar{L}_{h}u_{\gamma_{0}} 
+2\mu u_{i}\bar{L}_{h}u_{i}+2\kappa u\bar{L}_{h} u
$$
\begin{equation}
                                           \label{9.22.5}
+2\bar{a}_{k}[(\delta_{k}u_{\gamma_{0}}^{-})^{2}
+\mu\sum_{i}(\delta_{k}u_{i} )^{2}
+\kappa(\delta_{k}u )^{2}].
\end{equation}
 
Furthermore,
$$
0\geq \delta^{T}_{\tau}[U_{\gamma_{0}\mu}
+\kappa u^{2} ]=
 2u_{\gamma_{0}}^{-}\delta^{T}_{\tau}u_{\gamma_{0}}^{-} 
+2\mu u_{i}\delta^{T}_{\tau}u_{i}+2\kappa u\delta^{T}_{\tau} u
$$
$$ 
+2\tau [(\delta^{T}_{\tau}u_{\gamma_{0}}^{-})^{2}
+\mu\sum_{i}(\delta^{T}_{\tau}u_{i} )^{2}+
\kappa(\delta^{T}_{\tau}u)^{2}]
$$
\begin{equation}
                                           \label{9.22.6}
\geq
-2u_{\gamma_{0}}^{-}\delta^{T}_{\tau}u_{\gamma_{0}} 
+2\mu u_{i}\delta^{T}_{\tau}u_{i}+2\kappa u\delta^{T}_{\tau} u.
\end{equation}
We  multiply \eqref{9.22.6} by $\bar{r}$,
 add the result to \eqref{9.22.5},  
and use that 
$\bar{r}\geq0$, $\bar{c}\geq-K_{3}$,
$\bar{a}_{k}\geq\delta$, $|u|\leq K_{1}$,
and 
  $$ 
u(\bar{r}\delta^{T}_{\tau}u+\bar{L}_{h}u)=-u\bar{g} =\bar{c}u^{2}
-u(\bar{g}+\bar{c}u)\geq-K_{3}K_{1}-u(\bar{g}+\bar{c}u),
$$
where $|\bar{g}+\bar{c}u|
\leq \omega W +K_{3}
$
(recall \eqref{9.23_1}). Then we obtain
$$
u(\bar{r}\delta^{T}_{\tau}u+\bar{L}_{h}u)
\geq -\omega K_{1} W -2K_{1}K_{3}, 
$$
$$
 u_{\gamma_{0}}^{-}[\bar{r}\delta^{T}_{\tau}+\bar{L}_{h}]u_{\gamma_{0}} 
-\mu u_{i}[\bar{r}\delta^{T}_{\tau}+\bar{L}_{h}]u_{i}  
\geq
\delta\mu M
+\kappa  W(\delta-K_{1}\omega )   
-2
\kappa K_{1}K_{3}  ,
$$
where  
$$
M:=\sum_{i,k}(\delta_{k}u_{i} )^{2}.
$$

Since  $C\geq4$ in 
Assumption \ref{assumption 10.1.1},  we have $K_{1}\omega\leq\delta/4$
which along with  \eqref{9.23_1} leads to
$$
 W(\delta-K_{1}\omega )-2 K_{1}K_{3} \geq
(3/4)\delta W-2K_{1}K_{3}\geq(\delta/2)W,
$$  
\begin{equation}   
                                           \label{9.22.4}
 u_{\gamma_{0}}^{-}[\bar{r}\delta^{T}_{\tau}+\bar{L}_{h}]u_{\gamma_{0}} 
-\mu u_{i}[\bar{r}\delta^{T}_{\tau}+\bar{L}_{h}]u_{i}  
\geq
\delta\mu M
+(\delta/2)\kappa  W.
\end{equation}

On the other hand, owing to \eqref{9.23_1}, \eqref{9.21.4},
and Lemma \ref{lemma 9.21.1}
\begin{equation}
                                           \label{9.22.3}
P_{\gamma_{0}\mu}[\bar{r}\delta_{\tau}^{T}u+\bar{L}_{h}u+\bar{g}]
\leq0.
\end{equation}
Here due to \eqref{9.22.4}
 $$
P_{\gamma_{0}\mu}[\bar{r}\delta_{\tau}^{T}u+\bar{L}_{h}u]
= u_{\gamma_{0}}^{-}[\bar{r}\delta^{T}_{\tau}+\bar{L}_{h}]u_{\gamma_{0}} 
-\mu u_{i}[\bar{r}\delta^{T}_{\tau}+\bar{L}_{h}]u_{i}+I
$$
$$
\geq \delta\mu M
+(1/2)\delta \kappa W +I ,
$$
where
$$
I:=[u_{\gamma_{0}}^{-}\gamma_{0i}-\mu u_{i}]
(\delta_{i}\bar{a}_{k})[
\Delta_{k}u +h\Delta_{k}u_{i}].
$$

Below by $C$ we denote   generic constants depending only on
$d_{1}$.
It follows from   the estimates
$$
|\delta_{i}\bar{a}_{k}|\leq K_{3}+\omega|u_{i}|,\quad
|h\Delta_{k}u_{i}|=|\delta_{k}u_{i}
+\delta_{-k}u_{i}|\leq2
 M^{1/2},
$$
$$
|\Delta_{k}u|\leq   M^{1/2},
\quad
|u_{\gamma_{0}}^{-}\gamma_{0i} -\mu u_{i}|\leq CW^{1/2},
$$
 that 
$$
|I|\leq C  W^{1/2} M^{1/2}(K_{3}+\omega W^{1/2})
\leq NW+(1/2)\delta\mu M+C\delta^{-1}\omega^{2}W^{2}.
$$
 Hence,  
\begin{equation}
                                           \label{9.22.7}
P_{\gamma_{0}\mu}[\bar{r}\delta_{\tau}^{T}u+\bar{L}_{h}u]
\geq (1/2)\delta\mu M
+(1/2)\delta \kappa W-NW-C\delta^{-1}\omega^{2}W^{2}.
\end{equation}

To estimate $ P_{\gamma_{0}\mu}\bar{g}  $
recall that $\bar{c}\geq-N$, 
$U_{\gamma_{0}\mu}\geq0$ and observe that
$$
-P_{\gamma_{0}\mu}(\bar{c}u)=-\bar{c}P_{\gamma_{0}\mu}u
-[u_{\gamma_{0}}^{-}\gamma_{0k}-\mu u_{k}](T_{h,\ell_{k}}u)
\delta_{k}\bar{c} 
$$
$$
=\bar{c} U_{\gamma_{0}\mu}-[u_{\gamma_{0}}^{-}\gamma_{0k}-\mu
u_{k}](T_{h,\ell_{k}}u)
\delta_{k}\bar{c}
 \geq -NW-NW^{1/2} 
\geq -NW ,
$$
where the last inequality follows from \eqref{9.23_1}.
Furthermore,
$$
h\delta_{i}(\bar{g}+\bar{c} u)=  \bar{f} (T_{h,\ell_{i}} u_{k} ,
T_{h,\ell_{i}} u ,t_{0},x_{0}+h\ell_{i})
-\bar{f} (u_{k},
u ,t_{0},x_{0} ).
$$
Owing to  Assumption \ref{assumption 9.27.4},
\eqref{9.23_1},
and  the mean value theorem (this is the place,
where one cannot assume that \eqref{5.9.8} holds
 only for large $p$)
$$
|\delta_{i}(\bar{g}+\bar{c} u)|\leq C M^{1/2}
(\omega W^{1/2} +K_{3})
$$
$$
+C \big[
W^{1/2}(\omega W  +K_{3})
+\omega W^{3/2}+K_{3}\big].
$$

Note that the coefficients of $\delta_{i}(\bar{g}+\bar{c} u)$
in $P_{\gamma_{0}\mu}(\bar{g}+\bar{c}u)$ are dominated by $W^{1/2}$ and
for any $\rho>0$
$$
\omega M^{1/2}W\leq \rho M+\rho^{-1}\omega^{2} W^{2},
\quad
K_{3}M^{1/2}W^{1/2}\leq\rho M+\rho^{-1}NW.
$$
Therefore,
$$
|P_{\gamma_{0}\mu}(\bar{g}+\bar{c}u)|\leq
(1/2)\delta\mu M+C[\delta^{-1}\omega^{2}+\omega] W^{2}+NW,
$$
$$
P_{\gamma_{0}\mu}\bar{g}\geq
-(1/2)\delta\mu M-C[\delta^{-1}\omega^{2}+\omega] W^{2}-NW
$$
and \eqref{9.22.3} and \eqref{9.22.7} yield
$$
(1/2)\delta\kappa W\leq C[\delta^{-1}\omega^{2}+\omega] W^{2}+NW,
$$
\begin{equation}
                                                          \label{10.2.1}
2\delta\kappa  \leq C_{1}[\delta^{-1}\omega^{2}+\omega] W +N_{1}. 
\end{equation}
Since $W\leq\mu^{-1}4K_{1}^{2}\kappa$ (see \eqref{9.21.1}),
we see that if $C$ in Assumption \ref{assumption 10.1.1} is such that
$$
C\geq 4\mu^{-1}C_{1},
$$
then (recall that $C\geq4$ and $\omega K_{1}\leq\delta$)
$$
C_{1}[\delta^{-1}\omega^{2}+\omega] W
\leq C K_{1}^{2}[\delta^{-1}\omega^{2}+\omega]\kappa
\leq C[K_{1}\omega+K_{1}^{2}\omega]\kappa\leq\delta\kappa.
$$
In this case \eqref{10.2.1} allows us to conclude that
$\kappa\leq\delta^{-1}N_{1}$
and we get the assertion of the theorem from \eqref{9.21.1}.
The theorem is proved.

\mysection{Conditional estimates of the second-order differences}

                                                  \label{section 10.7.1}
                                                 
In this section we suppose that all the assumptions of Theorem 
\ref{theorem 10.4.1}
are satisfied apart from Assumptions \ref{assumption 10.3.2}
and \ref{assumption 10.7.1}.
The notation in this section are somewhat different from
Sections \ref{section 10.3.1} and \ref{section 10.20.1}.
Of course, we use our
 basic notation from Section \ref{section 5.9.2}, for instance,
$\xi$ and $\lambda$ are defined in \eqref{5.9.7}.

For $\varepsilon\in(0,1]$ set (observe that now $\gamma_{k}=\gamma_{-k}$)
$$
\Gamma(\varepsilon)=\{\gamma=(\gamma_{k}:k=\pm1,...,\pm d_{1}):
\gamma_{k}=\gamma_{-k},\varepsilon\leq\gamma_{k}\leq\varepsilon^{-1}
,\forall k\}.
$$
Fix a constant $\mu=\mu(d_{1},\varepsilon)>0$ such that 
$$
16d_{1}^{2}\mu\leq\varepsilon^{2}.
$$
In this section the indices $i,j,k,p,q$ run through
$\{\pm1,...,\pm d_{1}\}$.  
The main result of this section is the following.
 
\begin{theorem}
                                                 \label{theorem 10.3.1}  

Assume that
\begin{equation}
                                                 \label{10.6.8} 
3\mu\max_{Q}\sum_{i,j}\big[
(\xi\delta_{h,\ell_{j}}\delta_{h,\ell_{i}}u)^{-}\big]^{2}
\leq
 \max_{\Gamma(\varepsilon)\times Q} 
\big[(\xi\sum_{i}\gamma_{i}\Delta_{h,\ell_{i}}u)^{-}\big]^{2} .
\end{equation}

Then there exists a constant  $N=N(\varepsilon,\mu,d_{1},K_{0})$
 such that
if $ \lambda \geq N$, then in $Q|_{0}$
for $i,j=\pm 1,...,\pm d_{1} $ we have
$$
|\delta_{h,\ell_{j}}\delta_{h,\ell_{i}}u|
\leq 
N^{*}e^{m^{+}(T+\tau)}(1+\max_{\partial_{2}Q,p,q}
|\xi_{(-)}\delta_{h,\ell_{p}}\delta_{h,\ell_{q}}u|
$$
\begin{equation}
                                                 \label{8_27.3}
+
\max_{\bar{Q}}|\xi_{(-)}u|
+\max_{\bar{Q},p}|\xi_{(-)}\delta_{h,\ell_{p}}u|),
\end{equation}
where  $N^{*}=N^{*}(h_{0},\varepsilon,\mu,d_{1},K_{0},K_{3})$.
\end{theorem}

Below in the section by $N$ and $N^{*}$ we denote
generic constants of the same type as in the
theorem. As before,
we use the abbreviated notation
$$
\Delta_{i}=\Delta_{h,\ell_{i}},\quad \delta_{i}=\delta_{h,\ell_{i}}.
$$

Introduce     $v=\xi u$
as in Section \ref{section 10.3.1}  
and fix a constant  
  $\nu\geq1$.  Set
$$
P_{\gamma}=\gamma_{i}\Delta_{i},
\quad v_{\gamma}=P_{\gamma}v,\quad v_{i}=\delta_{i}v,\quad
v_{ij}=\delta_{j}\delta_{i}v,
$$  
$$
P_{\gamma\mu\nu}\phi=v^{-}_{\gamma}P_{\gamma} 
\phi+\mu
v_{ij}^{-}\delta_{j}\delta_{i}\phi-\nu v_{i}
\delta_{i}\phi,
$$
$$
W_{1}=\sum_{i}v_{i}^{2},\quad
W_{2}=\sum_{i,j}[v_{ij}^{-}]^{2},\quad V_{\gamma\mu\nu}=[v_{\gamma}^{-}]^{2}
+\mu W_{2}+\nu W_{1}.
$$

Observe that this time again $P_{\gamma\mu\nu}v=-V_{\gamma\mu\nu}$
and also note that \eqref{10.6.8} is equivalent to the following

\begin{equation}
                                                 \label{8_27.1}
3\mu\max_{Q}W_{2}\leq
 \max_{\Gamma(\varepsilon)\times Q} (v_{\gamma}^{-} )^{2} .
\end{equation}

We 
introduce $(\gamma_{0},t_{0},x_{0})$ as a point
in $\Gamma(\varepsilon)\times \bar{Q}$ maximizing
$V_{\gamma\mu\nu}$ and first prove  few auxiliary results.
Below, as usual, we drop the arguments $(t_{0},x_{0})$.
\begin{lemma}
                                           \label{lemma 8.25.3}
(i) For $(t,x)\in Q^{o}_{1}$ and any $i,j$  
\begin{equation}
                                                 \label{8_27.2}
|v_{ij}(t,x)|\leq \max_{Q}W^{1/2}_{2}.
\end{equation}

(ii)
If \eqref{8_27.1} holds and $(t_{0},x_{0})\in Q$  and
\begin{equation}
                                                 \label{10.5.1}
\nu \max_{Q}W_{1}\leq \mu \max_{Q}W_{2}  ,
\end{equation}
then at $(t_{0},x_{0})$
\begin{equation}
                                                 \label{8.25.6}
 \mu \max_{Q}W_{2} \leq
 [v_{\gamma_{0}}^{-}]^{2},\quad
 \nu \max_{Q}W_{1} \leq
 [v_{\gamma_{0}}^{-}]^{2}.
\end{equation}
Furthermore, if additionally, 
$$
h\sqrt{\nu}\leq\varepsilon,
$$  
 then
the operator $P_{\gamma_{0}\mu\nu}$
respects the maximum principle at $(t_{0},x_{0})$
relative to $\Lambda_{0}+\Lambda_{0}$, that is,
for any   function $\phi$ such that
 $\phi(x_{0})\geq\phi(x_{0}+\eta)$ for all
$\eta\in\Lambda_{0}+\Lambda_{0}$, we have $P_{\gamma_{0}\mu\nu}\phi(t_{0},x_{0})
\leq0$.

\end{lemma}

Proof. (i) Obviously $v_{-i,j}^{-}\leq W_{2}^{1/2}$ on $Q$.
Since (no summation in $i$) $T_{h,\ell_{i}}v_{-i,j}^{-}=v_{i j}^{+}$,
we get $v_{ij}^{+}(t,x)\leq W_{2}^{1/2} (t,x+h\ell_{i}) $. This proves
(i).

(ii) The second estimate in \eqref{8.25.6} follows from the first
one and \eqref{10.5.1}. Assuming that
the first estimate
in \eqref{8.25.6} does not hold, we obtain  at $(t_{0},x_{0})$
$$
V_{\gamma_{0}\mu\nu} < 2\mu\max_{Q} W_{2}+\nu\max_{Q}
W_{1}\leq 3\mu\max_{Q} W_{2} ,
$$
$$
\max_{\Gamma(\varepsilon)\times
Q}[v_{\gamma }^{-}]^{2}\leq V
_{\gamma_{0}\mu\nu} <3\mu\max_{Q}W_{2}
$$
contrary to \eqref{8_27.1}. This proves \eqref{8.25.6}.

 To prove the last assertion of the lemma we take a function
with described properties and without loss of generality assume
that $\phi(x_{0})=0$.
We also note that
$$
h^{2}\delta_{i}\delta_{j}\phi(x_{0})
=\phi(x_{0}+h\ell_{i}+h\ell_{j})-\phi(x_{0}+h\ell_{i})-
\phi(x_{0}+h\ell_{j})
$$
$$
+\phi(x_{0})\leq-\phi(x_{0}+h\ell_{i})-
\phi(x_{0}+h\ell_{j})
$$
and 
$\Delta_{i}\phi(x_{0})\leq0$. Therefore, as usual
dropping the arguments $(t_{0},x_{0})$ in  $v_{...}$, 
we infer from
\eqref{8.25.6} that
$$
h^{2}P_{\gamma_{0}\mu\nu}\phi(t_{0},x_{0})\leq v_{\gamma_{0}}^{-}
\varepsilon
\sum_{i}(\phi(x_{0}+h\ell_{i})
+\phi(x_{0}-h\ell_{i}))
$$
$$
-\mu v^{-}_{ij}(\phi(x_{0}+h\ell_{i})+
\phi(x_{0}+h\ell_{j}))-h\nu v_{i}\phi(x_{0}+h\ell_{i})
$$
$$
\leq 
v_{\gamma_{0}}^{-}\big[\varepsilon\sum_{i}(\phi(x_{0}+h\ell_{i})
+\phi(x_{0}-h\ell_{i}))
$$
$$
-\sqrt{\mu}\sum_{i,j}(\phi(x_{0}+h\ell_{i})+
\phi(x_{0}+h\ell_{j}))-h\sqrt{\nu}\sum_{i}\phi(x_{0}+h\ell_{i})\big]
$$
$$
=v_{\gamma_{0}}^{-}\sum_{i}\phi(x_{0}+h\ell_{i})[
2\varepsilon-2d_{1}\sqrt{\mu}-h\sqrt{\nu}].
$$
The last expression is less than zero 
in light of the fact that 
$4d_{1}\sqrt{\mu}\leq\varepsilon$ and
$h\sqrt{\nu}\leq\varepsilon$.
 The lemma is proved.
\begin{remark}
This lemma can be generalized to the case when $\ell_{k}$'s
come with different $h_{k}$'s, but the $h_{k}$'s
should be comparable. This is the reason why
in Theorem \ref{theorem 10.4.1} we do not include
$h_{d_{1}+1}$ and $\ell_{d_{1}+1}$.
\end{remark}
Set
$$
Z_{k}:=v_{\gamma_{0}}^{-}\Delta_{k}v_{\gamma_{0}}
+ \mu v_{ij}^{-}\Delta_{k}v_{ij}- \nu v_{i}\Delta_{k}
v_{i},
$$
$$
z_{k}:=v_{\gamma_{0}}^{-}\delta_{k}v_{\gamma_{0}}
+\mu v_{ij}^{-}\delta_{k}v_{ij}
-\nu v_{i}v_{ki},
$$
$$
R^{\gamma}_{k}:=\big[\delta_{k}v_{\gamma_{0}}^{-}\big]^{2}
,\quad R^{\mu}_{k}:=\sum_{i,j}\big[\delta_{k}v_{ij}^{-}\big]^{2}
,\quad R^{\nu}_{k}:=\sum_{i}v_{ki}^{2}.
$$
These objects evaluated at $(t_{0},x_{0})$ will be extensively
used below in the section.
 \begin{lemma}
                                                  \label{lemma 8.17.1}
If $(t_{0},x_{0})\in Q^{o}_{1}$, then
at $(t_{0},x_{0})$ we have for any $k$
$$
2Z_{k}\geq-2v_{\gamma_{0}}^{-}\Delta_{k}v_{\gamma_{0}}^{-}
-2\mu v_{ij}^{-}\Delta_{k}v_{ij}^{-}
-2\nu v_{i}\Delta_{k}
v_{i}
$$
\begin{equation}
                                                     \label{8.17.3}
\geq R^{\gamma}_{k}+R^{\gamma}_{-k}+\mu R^{\mu}_{k}+
\mu R^{\gamma}_{-k}
+\nu R^{\nu}_{k}+\nu R^{\nu}_{-k} \geq0.  
\end{equation}
 
Furthermore, for any $\alpha\in A$ and $h\leq h_{0}/2$
\begin{equation}
                                                     \label{10.13.1}
4(a^{\alpha}_{k}Z_{k}+
b^{\alpha}_{k}z_{k})\geq 
2a^{\alpha}_{k}Z_{k}+a^{\alpha}_{k}(R^{\gamma}_{k}
+\mu R^{\mu}_{k}+\nu R^{\nu}_{k}).
\end{equation}

\end{lemma}

Proof. 
The first inequality in \eqref{8.17.3} 
follows from Lemma \ref{lemma 6.18.2}. To prove the second one
it suffices to observe that
$$
0\geq\Delta_{k}V_{\gamma_{0}\mu\nu}
=2v_{\gamma_{0}}^{-}\Delta_{k}v_{\gamma_{0}}^{-}
+2\mu v_{ij}^{-}\Delta_{k}v_{ij}^{-}
+2\nu v_{i}\Delta_{k} v_{ i} 
$$
$$
+ R^{\gamma}_{k}+R^{\gamma}_{-k}+\mu R^{\mu}_{k}+
\mu R^{\mu}_{-k}
+\nu R^{\nu}_{k}+\nu R^{\nu}_{-k}.
$$

Next, using that 
 $h\leq h_{0}/2$, 
by \eqref{8.17.3} 
and  Assumption \ref{assumption 9.27.3} we get
$$
0\geq (a^{\alpha}_{k}\Delta_{k}
+2b^{\alpha}_{k}\delta_{k})V_{\gamma_{0}\mu\nu}
=2v_{\gamma_{0}}^{-}[a^{\alpha}_{k} \Delta_{k}
+2b^{\alpha}_{k}\delta_{k}]v_{\gamma_{0}}^{-}
$$
$$
+ 2\mu v_{ij}^{-}[a^{\alpha}_{k} \Delta_{k}
+2b^{\alpha}_{k}\delta_{k}]v_{ij}^{-}
+ 2\nu v_{i} [a^{\alpha}_{k} \Delta_{k}
+2b^{\alpha}_{k}\delta_{k}]v_{ i}  
$$
$$
+2a^{\alpha}_{k}[
R^{\gamma}_{k} +\mu R^{\mu}_{k} 
+\nu R^{\nu}_{k} ] +
2hb^{\alpha}_{k}[
R^{\gamma}_{k} +\mu R^{\mu}_{k} 
+\nu R^{\nu}_{k} ]
$$
$$
\geq-2v_{\gamma_{0}}^{-}[a^{\alpha}_{k} \Delta_{k}
+2b^{\alpha}_{k}\delta_{k}]v_{\gamma_{0}} 
$$
$$
- 2\mu v_{ij}^{-}[a^{\alpha}_{k} \Delta_{k}
+2b^{\alpha}_{k}\delta_{k}]v_{ij} 
+ 2\nu v_{i} [a^{\alpha}_{k} \Delta_{k}
+2b^{\alpha}_{k}\delta_{k}]v_{ i}  
$$
$$
+2a^{\alpha}_{k}[
R^{\gamma}_{k} +\mu R^{\mu}_{k} 
+\nu R^{\nu}_{k} ] +
2hb^{\alpha}_{k}[
R^{\gamma}_{k} +\mu R^{\mu}_{k} 
+\nu R^{\nu}_{k} ]
$$
$$
\geq-2a^{\alpha}_{k}Z_{k}-4b^{\alpha}_{k}z_{k}
+a^{\alpha}_{k}[
R^{\gamma}_{k} +\mu R^{\mu}_{k} 
+\nu R^{\nu}_{k} ]
$$
 and \eqref{10.13.1} follows.  
  The lemma is proved.

In the following lemma we do the most
important step in the proof of Theorem \ref{theorem 10.3.1}. Set
$$
\bar{W}_{1}:=\max_{Q}W_{1}.
$$
\begin{lemma}
                                                 \label{lemma 10.6.1}
Under the assumptions of Theorem \ref{theorem 10.3.1}
  there are constants $N$, $N^{*}$,
 $$
\nu=\nu ^{*}(h_{0},\varepsilon,\mu,d_{1},K_{0},K_{3})\geq1,
\quad
h^{*}=h^{*}(h_{0},\varepsilon,\mu,d_{1},K_{0},K_{3})>0
$$
such that   $h^{*}\leq h_{0}/2$ and,
if condition \eqref{10.5.1}
is satisfied and $h\in(0,h^{*}]$ and $(t_{0},x_{0})\in Q^{o}_{2}$,
then at $(t_{0},x_{0})$ for any $\alpha\in A$  
we have
\begin{equation}
                                                    \label{10.6.2}
J:=P_{\gamma_{0}\mu\nu}(a^{\alpha}_{k}\Delta_{k}
v +b^{\alpha}_{k}\delta_{k}v)\geq-N[v^{-}_{\gamma_{0}}]^{2}
-N^{*}\bar{W}_{1} .
\end{equation}
 
\end{lemma} 

Proof. We fix an $\alpha\in A$ and drop the superscript $\alpha$
for convenience.
 By \eqref{10.24.1}  (no summation in $k,i,j$)
$$
\delta_{j}\delta_{i}(b _{k}v_{k})=
b _{k}\delta_{k}v_{ij}+
(\delta_{j}b _{k})v_{ki}+
(\delta_{i}b _{k})v_{kj}
$$
$$
+h[(\delta_{i} +\delta_{j})
b _{k}]
 \delta_{k}v_{ij}+(\delta_{j}\delta_{i}b _{k})
T_{h,\ell_{i}+\ell_{j}}v_{k}.
$$
Also by using \eqref{10.24.1} and the formulas
$ a _{k}= a _{-k}$
and $h\Delta_{k}=\delta_{k}+\delta_{ -k }$
and summing with respect to $k$ (but not in $i,j$)
we get
$$
\delta_{j}\delta_{i}
( a _{k}\Delta_{k}v)
= a _{k}\Delta_{k}v_{ij}+
(\delta_{j} a _{k})\Delta_{ k }v_{ i}
+
(\delta_{i} a _{k})\Delta_{k}v_{j}
$$
$$
+2\big[(\delta_{j}+\delta_{i}) a _{k}\big]
\delta_{k}v_{ij}
+(\delta_{j}\delta_{i}
 a _{k})T_{h,\ell_{j}+\ell_{i}}\Delta_{k}v . 
$$
While applying this formula to $\Delta_{i}$
it is also useful to observe
that
$$
\big[(\delta_{ -i }+\delta_{i}) a _{k}\big]
\delta_{k}v_{i,-i}
=-h(\Delta_{i} a _{k})\Delta_{i}\delta_{k}v
=-(\Delta_{i} a _{k})(v_{ki}+v_{k,-i}).
$$
Hence, (recall that $\gamma_{0i}=\gamma_{0,-i}$ and $\ell_{i}=-\ell_{-i}$)
$$
P_{\gamma_{0}}
(
a _{k}\Delta_{k}v+b _{k}\delta_{k}v)
=-\gamma_{0i}\delta_{ -i }\delta_{i}
( a _{k}\Delta_{k}v
+b _{k}\delta_{k}v)
$$
$$
= (a _{k}\Delta_{k}+b _{k}\delta_{k})v_{\gamma_{0}} 
-2\gamma_{0i}
(\delta_{ -i } a _{k})\Delta_{ k }v_{i}
+4\gamma_{0i}(\Delta_{i} a _{k})v_{ki}
$$
$$
+\gamma_{0i}(\Delta_{i}
 a _{k}) \Delta_{k}v
+2\gamma_{0i}(\delta_{i}b _{k})v_{ki}
+\gamma_{0i}(\Delta_{i}b _{k})v_{k}.
$$
Also everywhere
$$
v_{ij}^{-}
 \delta_{j}\delta_{i}( a _{k}\Delta_{k}
+b _{k}\delta_{k})v
=v_{ij}^{-}( a _{k}\Delta_{k}
+b _{k}\delta_{k})v_{ij}
+2v_{ij}^{-}(\delta_{j} a _{k})
\Delta_{k}v_{i}
$$
$$
+hv_{ij}^{-}\big[(\delta_{j}+\delta_{i}) a _{k}\big]
\Delta_{k}v_{ij}
+
v_{ij}^{-}(\delta_{j}\delta_{i}
 a _{k})T_{h,\ell_{j}+\ell_{i}}\Delta_{k}v 
$$
$$
+2v_{ij}^{-}(\delta_{j}b _{k})v_{ki}
+2hv_{ij}^{-}(\delta_{i}b _{k})\delta_{k}v_{ij} 
+v_{ij}^{-}(\delta_{j}\delta_{i}b _{k})
T_{h,\ell_{i}+\ell_{j}}v_{k}.
$$

Therefore at $(t_{0},x_{0})$ we have  
$$
J=  a _{k}Z_{k}+b _{k}z_{k}+I_{1}+...+I_{4} ,
$$
 where 
$$
I_{1}=-2\gamma_{0i}v_{\gamma_{0}}^{-}
(\delta_{-i} a _{k})\Delta_{  k }v_{i}
+2\mu v_{ij}^{-}(\delta_{j} a _{k})
\Delta_{k}v_{i},
$$
$$
I_{2}=\mu hv_{ij}^{-}\big[(\delta_{j}+
\delta_{i}) a _{k}\big]
\Delta_{k}v_{ij},
$$ 
$$
I_{3}=v_{\gamma_{0}}^{-}[4\gamma_{0i}(\Delta_{i} a _{k})v_{ki}+
\gamma_{0i}(\Delta_{i} a _{k})\Delta_{k}v]
$$
$$
+\mu
v_{ij}^{-}(\delta_{j}\delta_{i}
 a _{k})T_{h,\ell_{j}+\ell_{i}}\Delta_{k}v,
$$
$$
I_{4}= 2v_{\gamma_{0}}^{-}\gamma_{0i}(\delta_{i}b _{k})v_{ki}
+v_{\gamma_{0}}^{-}\gamma_{0i}(\Delta_{i}b _{k})v_{k}
$$
$$
+2\mu v_{ij}^{-}(\delta_{j}b _{k})v_{ki}
+2\mu hv_{ij}^{-}(\delta_{i}b _{k})\delta_{k}v_{ij} 
+\mu v_{ij}^{-}(\delta_{j}\delta_{i}b _{k})
T_{h,\ell_{i}+\ell_{j}}v_{k}
$$
$$
-\nu v_{i}(\delta_{i} a _{k})
(\Delta_{k}v+2v_{ki})-
\nu v_{i}(\delta_{i}b _{k})T_{h,\ell_{i}}v_{k}.
$$

For $h\leq h_{0}/2$  
it follows by Lemma \ref{lemma 8.17.1} that
\begin{equation}
                                                        \label{8_27.7}
4J\geq 2 a _{k}
Z_{k}+
 a _{k}(R^{\gamma}_{k}
+\mu R^{\mu}_{k}+\nu R^{\nu}_{k})
+4I_{1}+...+4I_{4} .
\end{equation}

{\em Estimating $I_{1}$.\/} 
Note that  owing to  \eqref{8.25.6} 
$$
|4I_{1}|
\leq N v_{\gamma_{0}}^{-}(\sqrt{ a _{k}}+h)
\sum_{i}|\Delta_{k}v_{i}|
$$
and by \eqref{6.18.6}  
$$
Nv_{\gamma_{0}}^{-} \sqrt{ a _{k}} 
\sum_{i}|\Delta_{k}v_{i}|
\leq N v_{\gamma_{0}}^{-}\sqrt{ a _{k}}
\sum_{i}\big[|\delta_{-k}v_{ki}^{-}|+  
|\delta_{k}v_{-k,i}^{-}|\big]
$$
$$
\leq N(v_{\gamma_{0}}^{-})^{2}+
(1/3)\mu a _{k}R^{\mu}_{k}.
$$
Furthermore,
by the formula $h\Delta_{k}
=\delta_{k}+\delta_{h,\ell_{-k}}$
and Lemma \ref{lemma 8.25.3}
we obtain
$$
v_{\gamma_{0}}^{-}h
\sum_{k,i}|\Delta_{k}v_{i}|
\leq
2v_{\gamma_{0}}^{-} 
\sum_{k,i}|v_{ki}|
\leq Nv_{\gamma_{0}}^{-}\max_{Q}W^{1/2}_{2}
\leq N(v_{\gamma_{0}}^{-})^{2}.
$$
Thus,
\begin{equation}
                                                        \label{8_27.4}
|4I_{1}| \leq N(v_{\gamma_{0}}^{-})^{2}+
(1/3)\mu a _{k}
 R^{\mu}_{k}.
\end{equation}

{\em Estimating $I_{2}$.\/} Observe that
$$
|4I_{2}|\leq Nhv_{ij}^{-}(\sqrt{ a _{k}}+h)
|\Delta_{k}v_{ij}|\leq  I_{21}+I_{22} ,
$$
where (see Lemma \ref{lemma 8.25.3} and recall that $(t_{0},x_{0})\in
Q^{0}_{2}$)  
$$
I_{21}=N\sum_{k}h^{2}v_{ij}^{-}|\Delta_{k}v_{ij}|
=N\sum_{k} v_{ij}^{-}|(T_{h,\ell_{k}}-2+T_{h,\ell_{-k}})v_{ij}|
\leq N\ (v_{\gamma_{0}}^{-})^{2} ,
$$
and by the formula $|\theta|=\theta+2\theta^{-}$,
$$
I_{22} =Nhv_{ij}^{-} \sqrt{ a _{k}} 
|\Delta_{k}v_{ij}|= Nhv_{ij}^{-}\sqrt{ a _{k}}
\Delta_{k}v_{ij}
$$
$$
+2N
 h  v_{ij}^{-}\sqrt{ a _{k}}  
(\Delta_{k}v_{ij})^{-}=Nh\sqrt{a_{k}}Z_{k}
+N
 h  v_{ij}^{-}\sqrt{ a _{k}}  
(\Delta_{k}v_{ij})^{-}+I_{23},
$$
where
$$
 I_{23} =-Nh 
  \sqrt{ a _{k}}(v_{\gamma_{0}}^{-}\Delta_{k}v_{\gamma_{0}} 
-  \nu
v_{i}\Delta_{k}v_{i}) .
$$

Furthermore, by Lemma  \ref{lemma 6.18.2}   
$$
N h  v_{ij}^{-}\sqrt{ a _{k}}  
(\Delta_{k}v_{ij})^{-}
\leq N h\sqrt{ a _{k}} v_{ij}^{-} | 
 \Delta_{k}v_{ij} ^{-}|
=N\sqrt{ a _{k}} v_{ij}^{-} | 
( \delta_{k}+ \delta_{-k})v_{ij} ^{-}|,
$$
which is majorated by the right-hand side of \eqref{8_27.4}.

To estimate $I_{23}$ we use Lemma \ref{lemma 6.18.2} to get
$$
-Nh 
  \sqrt{ a _{k}} v_{\gamma_{0}}^{-}\Delta_{k}v_{\gamma_{0}}
=-N
  \sqrt{ a _{k}} v_{\gamma_{0}}^{-}\delta_{k}v_{\gamma_{0}}
\leq
N \sqrt{ a _{k}}v_{\gamma_{0}}^{-}
\delta_{k}v_{\gamma_{0}}^{-}
$$
\begin{equation}
                                                    \label{10.6.1}
\leq N(v_{\gamma_{0}}^{-})^{2}
+(1/3) a _{k}R^{\gamma}_{k}.
\end{equation}
Furthermore, by assumption \eqref{10.5.1}
$$
Nh\nu
v_{i} \sqrt{ a _{k}}\Delta_{k}v_{i}=N\nu v_{i} \sqrt{ a _{k}}
v_{ki}
\leq N(v_{\gamma_{0}}^{-})^{2}+(1/3)\nu
a _{k}R^{\nu}_{k}.
$$
It follows that 
$$
I_{23}\leq N(v_{\gamma_{0}}^{-})^{2} +
(1/3) a _{k}\big[\delta_{k}v_{\gamma_{0}}^{-}]^{2}
+(1/3)\nu a _{k}R^{\nu}_{k}.
$$

Hence   
\begin{equation}
                                                        \label{8_27.6}
|4I_{2}|\leq  N(v_{\gamma_{0}}^{-})^{2} 
+Nh \sqrt{ a _{k}}
Z_{k}
+(1/3) a _{k}(R^{\gamma}_{k}+\mu R^{\mu}_{k}+\nu R^{\nu}_{k}).
\end{equation}

{\em Estimating $I_{3}$.\/}
We use the following result of simple computations
$$
T_{h,\ell_{j}+\ell_{i}}\Delta_{k}v
=-v_{k,-k}+v_{kj}+v_{-k,j}+v_{ki}+v_{-k,i}+h^{2}
 \Delta_{k} v_{ij}.
$$
This shows new terms entering $I_{3}$. All of them
apart from the last one are 
similar to the ones which are written explicitly in
the definition of $I_{3}$ and we show
how to estimate only one of them. By Assumption
\ref{assumption 10.2.2}  and Lemma \ref{lemma 8.25.3} we have 
$$
|\mu v_{ij}^{-}(\delta_{j}\delta_{i}
 a _{k})v_{k,-k}|\leq N(v_{\gamma_{0}}^{-})^{2}
+N^{*}|v_{\gamma_{0}}^{-}|\sum_{k,i}\sqrt{ a _{k}}|v_{ki}|
\leq N(v_{\gamma_{0}}^{-})^{2}+N^{*}
 a _{k}R^{\nu}_{k}.
$$
To estimate the remaining term in $I_{3}$ we proceed
as in estimating $I_{2}$. We have
$$
 \mu h^{2}  v_{ij}^{-}(\delta_{j}\delta_{i}
 a _{k})\Delta_{k}v_{ij}  
$$
$$
\leq N h^{2} v_{ij}^{-}
(N+N^{*}\sqrt{ a _{k}} )|\Delta_{k}v_{ij}|=I_{31}+I_{32}.
$$
Here  by Lemma \ref{lemma 8.25.3} and because $(t_{0},x_{0}) 
\in Q^{o}_{2}$
$$
 I_{31}=N h^{2}  v_{ij}^{-}|\Delta_{k}v_{ij}|
=N  v_{ij}^{-}|(T_{h,\ell_{k}}-2+T_{h,\ell_{-k}})v_{ij}|
\leq N
(v_{\gamma_{0}}^{-})^{2}.
$$
Next,
$$
I_{32}=N^{*}h^{2}
 v_{ij}^{-}\sqrt{ a _{k}} |\Delta_{k}v_{ij}|
=I_{321}+I_{322},
$$
with
$$
I_{321}=N^{*}_{1} h^{2}
 \mu v_{ij}^{-}\sqrt{ a _{k}}  \Delta_{k}v_{ij},
$$
$$
I_{322}=N^{*}
h^{2}v_{ij}^{-}\sqrt{ a _{k}}(\Delta_{k}v_{ij})^{-}
\leq
N^{*}
h^{2}v_{ij}^{-}\sqrt{ a _{k}}|\Delta_{k}v_{ij}^{-}|
$$
$$
=N^{*}
h v_{ij}^{-}\sqrt{ a _{k}}|\delta_{k}v_{ij}^{-}
+\delta_{-k}v_{ij}^{-}|
\leq N^{*}\bar{W}_{1} +(1/3)\mu  a _{k} 
R^{\mu}_{k},
$$
where the last inequality is true
 since $h v_{ij}^{-}\leq 2\bar{W}_{1}^{1/2}$.
Also observe that
$$
I_{321}=N^{*}_{1} h^{2}
  \sqrt{ a _{k}}Z_{k}
-N^{*}_{1} h^{2}v_{\gamma_{0}}^{-}\sqrt{ a _{k}}\Delta_{k}v_{\gamma_{0}}
+N^{*}_{1} h^{2}\nu v_{i}\sqrt{ a _{k}}\Delta_{k}v_{i},
$$
where
$$
N^{*}_{1}\nu h^{2}v_{i}\sqrt{ a _{k}}\Delta_{k}v_{i}
=2N^{*}_{1}\nu h v_{i}\sqrt{ a _{k}} v_{ki}
\leq N^{*}\nu^{2}h^{2}\bar{W}_{1}+a_{k}R^{\nu}_{k}
$$
and,
according to \eqref{10.6.1} and the inequality
 $h |v_{\gamma_{0}}|\leq N\bar{W}_{1}^{1/2}$,
$$
-N^{*}_{1}h^{2}\sqrt{ a _{k}} v_{\gamma_{0}}^{-}
\Delta_{k}v_{\gamma_{0}}
\leq N^{*} \bar{W}_{1}+(1/3)
 a _{k}R^{\gamma}_{k}.
$$
Therefore,
$$
I_{321}\leq N^{*}_{1}h^{2}\sqrt{ a _{k}}
Z_{k}
+
 N^{*}(1+\nu^{2}h^{2})\bar{W}_{1}+
 (1/3)
 a _{k}R^{\gamma}_{k}+a_{k}R^{\nu}_{k}.
$$

We can now specify $h^{*}$: we take
 $$
 N^{*}_{1}h^{*}\leq1\quad\text{and}\quad h^{*}\leq h_{0}/2.
$$ 
 Then, for $h\leq h^{*}$,
$$
I_{32}\leq  h\sqrt{ a _{k}}
Z_{k}
+
 N^{*}\nu^{2}\bar{W}_{1} 
+(1/3)\mu
 a _{k} 
R^{\mu}_{k}
+ (1/3)
 a _{k}R^{\gamma}_{k}+ a_{k}R^{\nu}_{k},
$$
$$
|4I_{3 }|\leq  h\sqrt{ a _{k}}
Z_{k}
+
 N^{*}\nu^{2}\bar{W}_{1}+
  N (v_{\gamma_{0}}^{-})^{2}
$$
\begin{equation}
                                                        \label{8_27.5}
+(1/3)
 a _{k}R^{\gamma}_{k}
+(1/3)\mu
 a _{k} 
R^{\mu}_{k}+N^{*}
 a _{k} R^{\nu}_{k}.
\end{equation}

{\em Estimating $I_{4}$.\/} By using  Lemma \ref{lemma 8.25.3}
we easily  see that
$$
|4I_{4}|\leq N(v_{\gamma_{0}}^{-})^{2}+N^{*}  \bar{W}_{1}+
\nu\bar{W}_{1}^{1/2} (N\sqrt{a _{k}}
+N^{*}h)\sum_{i}|v_{ki}|,
$$
 where
$$
N\nu\bar{W}_{1}^{1/2} \sqrt{a _{k}}
 \sum_{i}|v_{ki}| \leq
N\nu\bar{W}_{1}+(1/3)\nu a _{k}R^{\nu}_{k},
$$
$$
N^{*}\nu\bar{W}_{1}^{1/2}h \sum_{k,i}|v_{ki}|\leq
N^{*}\nu\bar{W}_{1}.
$$
It follows that
$$
|4I_{4}|\leq N(v_{\gamma_{0}}^{-})^{2}+ N^{*}\nu\bar{W}_{1}+
(1/3)\nu a _{k}R^{\nu}_{k}.
$$

 By combining this with \eqref{8_27.4}, \eqref{8_27.6}, and
\eqref{8_27.5}, recalling that $Z_{k}\geq0$,
and coming back to \eqref{8_27.7} we conclude 
$$
4J\geq  (2a _{k}-N_{1}h\sqrt{a _{k}})
Z_{k}
$$
\begin{equation}
                                                        \label{10.6.5}
 -N(v_{\gamma_{0}}^{-})^{2}
-N^{*}\nu^{2}\bar{W}_{1}+(\nu/3-N^{*}_{2}) 
a _{k}R^{\nu}_{k}.
\end{equation}

Now we 
specify $\nu=\nu^{*}$ by setting 
$$
\nu^{*}=1+3N^{*}_{2}
$$ 
and
finish the argument
 as in \cite{Kr06}. Namely, if $ 2a _{k}-N_{1}h\sqrt{ a _{k}}
\geq0$ for a $k$, then we can drop the term on the right
in \eqref{10.6.5} corresponding to this $k$ because 
$Z_{k}\geq0$. However, if
$ 2a _{k}-N_{1}h\sqrt{ a _{k}}\leq0$, then
$\sqrt{ a _{k}}\leq Nh $ and
$| a _{k}-N_{1}h\sqrt{ a _{k}}|\leq Nh^{2}$,
whereas
$$
 h^{2}Z_{k}=
h^{2} ( v_{\gamma_{0}}^{-}\Delta_{k} v_{\gamma_{0}}+
\mu v_{ij}^{-}\Delta_{k}v_{ij} -\nu v_{i}
\Delta_{k}v_{ki})
$$
$$
\leq
N\max_{Q}W_{2}+N\nu\bar{W}_{1}^{1/2}\max_{Q}W_{2}^{1/2}
$$
$$
\leq N\max_{Q}W_{2}+N\nu^{2}\bar{W}_{1}
\leq N(v^{-}_{\gamma_{0}})^{2}+N^{*}\bar{W}_{1}.
$$

This and \eqref{10.6.5} yield \eqref{10.6.2} and the lemma is proved.

{\bf Proof of Theorem \ref{theorem 10.3.1}}.  
Fix a constant $\nu$ according to Lemma \ref{lemma 10.6.1}
and first assume that
$(t_{0},x_{0})\in\partial_{2}Q$. Then
$$
\sqrt{\mu} \max_{\bar{Q}}W_{2}^{1/2}
\leq V_{\gamma_{0}\mu\nu}^{1/2}(t_{0},x_{0})
\leq N^{*}(\max_{\partial_{2}Q,i,j}|v_{ij}|+
\max_{\partial_{2}Q,i}|v_{i}|),
$$
which by   \eqref{8_27.2} yields similar estimate for
$$
\max_{Q^{o}_{1},i,j}|v_{ij}|.
$$
After that \eqref{8_27.3} is immediate (cf.~the end of the proof
of Theorem \ref{theorem 10.11.1}).

Therefore, in the rest of the proof we assume that
$$
(t_{0},x_{0})\in Q^{o}_{2}.
$$
Similarly, if \eqref{10.5.1} is violated,
there is nothing to prove. Hence, we may assume that
\eqref{10.5.1} holds. Finally,
we may assume that $h\leq h^{*}$,
where $h^{*}$ is taken from Lemma \ref{lemma 10.6.1}
and further reduced it if needed so as to satisfy
$h^{*}\sqrt{\nu}\leq\varepsilon$. Indeed, if $h\geq h^{*}$,
then in $Q^{o}_{1}$
$$
|v_{ij}|\leq 2(h^{*})^{-1}\max_{Q}|v_{i}|.
$$

After justifying these additional assumptions which allow us to use
the assertions of Lemmas \ref{lemma 8.25.3}
 and \ref{lemma 10.6.1} as long as $h\leq h^{*}$, we construct 
functions $\bar{r}(t)$, 
$\bar{a}_{k}
(t,x)$, $\bar{b}_{k}(t,x)$, $\bar{c}(t,x)$, $\bar{f}(p,\psi,t,x)
=\bar{f}(t,x)$
as in Section \ref{section 10.3.1} to get \eqref{9.7.01}
 and \eqref{9.7.02} satisfied.
Then, since $(t_{0},x_{0})\in Q^{o}_{2}$ and
\eqref{8_27.1} and \eqref{10.5.1} 
are valid and $h\sqrt{\nu}\leq\varepsilon$, by Lemma \ref{lemma 8.25.3}
at
$(t_{0},x_{0})$ we obtain
$$
P_{\gamma_{0}\mu\nu}
\big(e^{-m \tau}\bar{r}\delta^{T}_{\tau}v
 +\bar{a}_{k}\Delta_{k}v+\bar{b}_{k}
\delta_{k}v-(\bar{c}+\bar{r}c_{m})v+\xi\bar{f}\big)\leq0.
$$
 
The fact that $\bar{r}$, $\bar{a}_{k}$,
 and $\bar{b}_{k}$ are  limits of some $r^{\alpha}$,
$a_{k}^{\alpha}$, and $b^{\alpha}_{k}$,
allows us to assert that Lemma \ref{lemma 10.6.1} holds
with $\bar{r}$,
$\bar{a}_{k}$,  and $\bar{b}_{k}$ in place of $r^{\alpha}$,
 $a_{k}^{\alpha}$,
and $b^{\alpha}_{k}$, respectively.
Therefore,  
\begin{equation}
                                                        \label{10.6.6}
 -N(v_{\gamma_{0}}^{-})^{2}+
P_{\gamma_{0}\mu\nu}
\big(e^{-m \tau}\bar{r}\delta^{T}_{\tau}v
 -(\bar{c}+\bar{r}c_{m})v+\xi\bar{f}\big)
\leq N^{*}
\bar W_{1}  .
\end{equation} 

Here $P_{\gamma_{0}\mu\nu}
 \delta^{T}_{\tau}v\geq0$
 as right after \eqref{9.13.2}.
 Furthermore,
$$
P_{\gamma_{0}\mu\nu}((\bar{c}+\bar{r}c_{m})v)
=-(\bar{c}+\bar{r}c_{m})V_{\gamma_{0}\mu\nu}
+I_{1},
$$
where $I_{1}$ is a linear combination of products of two types:

(i) $v_{\gamma_{0}}^{-}$ or $v_{ij}^{-}$
times a difference operator applied to $\bar{c}$
times either $v$ or a first-order difference
operator applied to $v$--  the second and third factors
may be taken at a point different from $(t_{0},x_{0})$,
but their coefficients in the linear combination
are dominated by a constant $N$;

(ii) $v_{i}$ times a difference operator applied to $\bar{c}$
times either $v$ or a first-order difference
operator applied to $v$-- these terms
may be taken at a point different from $(t_{0},x_{0})$, the
coefficients of these terms are dominated by a constant
$N^{*}$ (recall that $\nu$ is entering $P_{\gamma\mu\nu}$).

Owing to Lemma \ref{lemma 8.25.3},
the absolute value  of the linear combination of the products of type (i)  
is less than
$$
Nv_{\gamma_{0}}^{-}K_{3}(\max_{Q}|v|+\bar{W}^{1/2}_{1})
\leq (v_{\gamma_{0}}^{-})^{2}
+N^{*}(\max_{Q}|v|^{2}+\bar W _{1}).
$$
The absolute value  of the 
the linear combination of the products of type
(ii)  is clearly less than
$$
N^{*}\bar{W}^{1/2}_{1}(\max_{Q}|v|+ \bar{W}^{1/2}_{1})
\leq N^{*}(\max_{Q}|v|^{2}+\bar W _{1}).
$$

Now from the above estimates, \eqref{10.6.6},   
  and the fact that $V_{\gamma_{0}\mu\nu}
\geq (v_{\gamma_{0}}^{-})^{2}$ and  $\bar{c}+\bar{r}c_{m}\geq\lambda$
 we conclude
\begin{equation}
                                                        \label{10.6.7}
 ( \lambda -N)(v^{-}_{\gamma_{0}})^{2}
+P_{\gamma_{0}\mu\nu}
\big(  \xi\bar{f}\big)\leq
 N^{*}(\max_{Q}|v|^{2}+\bar W _{1}).
\end{equation}

Finally, obviously
$$
P_{\gamma_{0}\mu\nu}(\xi\bar{f})\geq-N^{*}\xi
(\max_{i,j}|v_{ij}|+\nu\max_{i}|v_{i}|)
\geq-N^{*}\xi^{2}-(v^{-}_{\gamma_{0}})^{2}  
-  \bar{W}_{1}
$$
and we infer from \eqref{10.6.7} that
$$
 ( \lambda -N_{1})(v^{-}_{\gamma_{0}})^{2}
 \leq
 N^{*}(\max_{Q}|v|^{2}+\bar{W} _{1})+N^{*}\xi^{2}.
$$
We set the constant $N$ in  the statement of the theorem
to be $N_{1}+1$ and use Lemma \ref{lemma 8.25.3}
to conclude that in $Q^{o}_{1}\cap Q|_{0}$ for any $i,j$
$$
 |\delta_{j}\delta_{i}u |=
|v_{ij} |\leq Nv_{\gamma_{0}}^{-}(t_{0},x_{0})
\leq
 N^{*}(\max_{Q}|v| +\bar{W}_{1}^{1/2})+N^{*}\xi(t_{0}).
$$
This implies \eqref{8_27.3}
in $Q^{o}_{1}\cap Q|_{0}$. On the remaining part of $Q|_{0}$
estimate \eqref{8_27.3} is obvious and the theorem is proved.

\mysection{Proof of Theorem \protect\ref{theorem 10.4.1}}
                                                 \label{section 10.20.2}

We start with three auxiliary results.
Everywhere in this section the assumptions
of Theorem \ref{theorem 10.4.1} are supposed to be satisfied.
Recall that the set  $\cL$ is introduced in
\eqref{5.10.3}.
\begin{lemma}
                                               \label{lemma 8.24.1}
For any function $\phi$ and $l_{1},l_{2}\in  \Lambda_{0}$ 
we have  
$$
|\delta_{h,l_{1}}\delta_{h,l_{2}}\phi(0)|\leq 4
 \max(|\Delta_{h,\ell_{k}}\phi(x)|:|k|\leq d_{0},x\in  
(\Lambda_{0}+\cL)\cup\{0\}  )
$$
\begin{equation}
                                                \label{8.24.2}
+4\max(|\Delta_{h,\ell_{k}}\phi(x)|:d_{0}<|k|\leq d_{1},
x\in \Lambda_{0}\cup\{0\}) .
\end{equation}

\end{lemma}

Proof. 
Obviously we may assume that $h=1$.
Next, observe that
$$
\delta_{1,l_{1}}\delta_{1,l_{2}}\phi(0)
=(1/2)[\Delta_{1,l_{2}}\phi(l_{1})+\Delta_{1,l_{1}}\phi(l_{2})]
-(1/2)\Delta_{1,l_{1}-l_{2}}\phi(0),
$$
$$
\delta_{1,l_{1}}\delta_{1,l_{1}}\phi(0)=\Delta_{h,l_{1}}\phi(l_{1}),
\quad \delta_{1,l_{1}}\delta_{1,-l_{1}}\phi(0)
=-\Delta_{h,l_{1}}\phi(0).
$$
It follows that if $l_{1},l_{2}\in \cL$, then
$$
|\delta_{h,l_{1}}\delta_{h,l_{2}}\phi(0)|\leq  
 \max(|\Delta_{h,\ell_{k}}\phi(x)|:|k|\leq d_{0},x\in  
\cL\cup\{0\}  )
$$
$$
+ \max(|\Delta_{h,\ell_{k}}\phi(0)|:d_{0}<|k|\leq d_{1}) .
$$
We substitute here
$\phi(y+\cdot)$ in place of $\phi$ 
and use that $\Lambda_{0}+\cL\supset\cL$ since
$d_{1}\geq2$ and $\cL=-\cL$. Then  we see
that, if 
$y\in \Lambda_{0} \cup\{0\}$ and $l_{1},l_{2}\in\cL$,
then
$$
|\delta_{1,l_{1}}\delta_{1,l_{2}}\phi(y)|
\leq  
 \max(|\Delta_{h,\ell_{k}}\phi(x)|:|k|\leq d_{0},x\in  
 (\Lambda_{0}+ \cL)\cup\{0\}  )
$$
\begin{equation}
                                                \label{8.24.3}
+ \max(|\Delta_{h,\ell_{k}}\phi(x)|:d_{0}<|k|\leq d_{1},
x\in\Lambda_{0}\cup\{0\}) .
\end{equation}

In case $l_{1}=\zeta_{1}+\zeta_{2},l_{2}=\eta_{1}+\eta_{2}$ with
$\zeta_{1},\zeta_{2},\eta_{1},\eta_{2}\in \cL $
and $\zeta_{1}\ne\zeta_{2}$, $\zeta_{1}\ne-\zeta_{2}$,
$\eta_{1}\ne\eta_{2}$, $\eta_{1}\ne-\eta_{2}$
either $\zeta_{1}\ne\eta_{1}$ and $\zeta_{1}\ne-\eta_{1}$
or $\zeta_{1}\ne\eta_{2}$ and $\zeta_{1}\ne-\eta_{2}$.
The second possibility reduces to the first one by interchanging
$\eta_{1}$ and $\eta_{2}$. If the first possibility realizes,
then
 we use the formula
$ 
\delta_{1,\eta+\zeta}  
=T_{1,\eta}\delta_{1,\zeta} +\delta_{1,\eta} 
$ 
to obtain
$$
\delta_{1,l_{1}}\delta_{1,l_{2}}\phi(0)=(T_{1,\zeta_{1}}\delta_{1,\zeta_{2}}
+\delta_{1,\zeta_{1}})(T_{1,\eta_{1}}\delta_{1,\eta_{2}}
+\delta_{1,\eta_{1}})\phi(0)
$$
$$
=\delta_{1,\zeta_{2}}\delta_{1,\eta_{2}}\phi(\zeta_{1}+\eta_{1})
+\delta_{1,\zeta_{2}}\delta_{1,\eta_{1}}\phi(\zeta_{1})
+\delta_{1,\zeta_{1}}\delta_{1,\eta_{2}}\phi(\eta_{1})
+\delta_{1,\zeta_{1}}\delta_{1,\eta_{1}}\phi(0).
$$
Here $\zeta_{1}+\eta_{1}\in\Lambda_{0},\zeta_{1},\eta_{1}\in
\Lambda_{0}$, and $0\in\Lambda_{0}\cup\{0\}$.
Therefore, we get \eqref{8.24.2} from \eqref{8.24.3}. 

The remaining case that $l_{1}\in\cL$ and $l_{2}=\eta_{1}+\eta_{2}$
with $\eta_{i}$ as above is taken care of by setting $\zeta_{1}=0$
in the above calculations.
The lemma is
proved.

Before stating the next lemma we remind the reader that
the index
$k$ takes values in
$\{\pm1,...,\pm d_{1}\}$.
\begin{lemma}  
                                               \label{lemma 8.25.1}
For any  values of the arguments and $s>0$ we have
$$
\sum_{k}q_{k}^{+}\leq \frac{2d_{1}}{\delta}
\big[s^{-1}F(s\phi,sq_{k},sp_{k},s\psi )+K_{0}
\sum_{k}q_{k}^{-}
$$
$$
+K_{3}\big(\sum_{k}|p_{k}|+ |\psi|+\phi^{-}+s^{-1}\big)\big].  
$$
\end{lemma}

Indeed,  the expression in the brackets obviously
is bigger than
$$
\sup_{\alpha\in A} a_{k}^{\alpha}(t,x)q_{k}^{+},
$$
which in turn is bigger than $\delta q_{n}^{+}$
for each particular $n=\pm1,...,\pm d_{1}$.
 
Below we use the notation $\Gamma(\varepsilon)$
and $P_{\gamma}$
from Section \ref{section 10.7.1}.
\begin{lemma}
                                               \label{lemma 8.25.2}
 Let  $\theta\in(0,\varepsilon^{-1})$
and $\varepsilon\in(0,1]$  
 be such that 
\begin{equation}
                                                 \label{8.25.4}
2d_{1}K_{0}\kappa/\delta\leq1/2,\quad
\kappa:=\varepsilon(\theta+\varepsilon)/(1-\theta\varepsilon).
\end{equation}
Let $w$, $\psi$, $p_{k}$ be  functions on $Q$ and assume that
\begin{equation}
                                                 \label{8.25.2}
\max_{\Gamma(\varepsilon)\times
Q} (\xi P_{\gamma}w)^{-}
\leq\theta\max_{Q}\sum_{k} |\xi\Delta_{h,\ell_{k}}w|.
\end{equation}
Then in $Q$  we have
$$
\sum_{k}
|\xi\Delta_{h,\ell_{k}}w|\leq\frac{4d_{1}}{\delta}(1+\kappa)I,
$$
where
$$
I=
\max_{Q}[\xi F ( \phi,
\Delta_{h,\ell_{k}}w,p_{k},\psi) 
+K_{3}\xi\big(\sum_{|k|\leq d_{1}}|p_{k}|+ |\psi|+
\phi^{-}+1\big)\big].
$$
 
\end{lemma}

Proof. Set $\Phi^{\pm}=\sum_{k}(\xi\Delta_{h,\ell_{k}}w)^{\pm}$
and observe that due to  \eqref{8.25.2} 
$$
\varepsilon\Phi^{+}-\varepsilon^{-1}\Phi^{-}
\geq-\theta\max_{Q}\Phi^{+}-\theta\max_{Q}\Phi^{-}
$$
in $Q$.
Hence $(\theta+\varepsilon)\max_{Q
}\Phi^{+}\geq(\varepsilon^{-1}-\theta)
\max_{Q }\Phi^{-}$, that is
\begin{equation}
                                                 \label{8.25.3}
\max_{Q }\Phi^{-}\leq
\kappa\max_{Q }\Phi^{+}.
\end{equation}
By \eqref{8.25.3} and
Lemma \ref{lemma 8.25.1} with $s^{-1}\phi,s^{-1}\Delta_{h,\ell_{k}}
w,s^{-1}p_{k}$, and $s^{-1}\psi$ in place of
$\phi,q_{k},p_{k}$, and $\psi$, respectively, and $s^{-1}=\xi$
we find  that in $Q$
$$
\Phi^{+}\leq\frac{2d_{1}}{\delta}
\big[I+K_{0}\kappa\max_{Q}\Phi^{+}\big].
$$

Upon taking the maximums over $Q $ of both parts and taking into
account
\eqref{8.25.4} we get that $\Phi^{+}\leq (4d_{1}/\delta)I$ in $Q $,
which along  with \eqref{8.25.3} yield the result.
The lemma is proved.  

{\bf Proof of Theorem \ref{theorem 10.4.1}}. 
Here $k,i,j$ run
through $\pm1,...,\pm d_{1}$.
It is easy to see that 
one can find $\varepsilon=\varepsilon(\delta,d_{1},K_{0})\in(0,1]$ and
$\mu=\mu(\delta,d_{1},K_{0})>0$ in such a way that the conditions:
$$
16d_{1}^{2}\mu\leq\varepsilon^{2},\quad
16d_{1}^{2}(3\mu)^{1/2}=:\theta<\varepsilon^{-1}
$$
and \eqref{8.25.4} are satisfied.
We choose and fix appropriate $\varepsilon$ and $\mu$.

If  
\begin{equation}
                                             \label{10.27.1}
\max_{\Gamma(\varepsilon)\times Q}(\xi P_{\gamma}u)^{-}
\leq \theta\max_{Q}  \sum_{ k}|\xi \Delta_{h,\ell_{k}}u|,
\end{equation}
then, by taking into account that $u$ satisfies \eqref{10.26.1}
in $Q $,
 from
Lemma \ref{lemma 8.25.2} we obtain that
$|\xi\Delta_{h,\ell_{k}}u|$ are bounded in $Q$
by the right-hand side  of  
 \eqref{10.8.1}.
By combining this with Lemma \ref{lemma 8.24.1} we conclude that  
\eqref{10.8.1} is true in $Q|_{0}\cap Q^{o}_{2}$ (notice that 
$(\Lambda_{0}+\cL)\cup\{0\} 
\subset\Lambda_{0}+\Lambda_{0}$).
Of course, \eqref{10.8.1} is obvious on $Q|_{0}\cap\partial_{2}Q$.
 
If  \eqref{10.6.8} 
holds, then we get  \eqref{10.8.1} from 
Theorem~\ref{theorem 10.3.1}.
In the remaining case  both \eqref{10.27.1}  and
\eqref{10.6.8} 
  are violated and
$$
\theta H: 
=\theta\max_{Q}  \sum_{ k}|\xi\Delta_{h,\ell_{k}}u|
$$
$$
\leq\max_{\Gamma(\varepsilon)\times Q}(\xi P_{\gamma}u)^{-}
\leq (3\mu)^{1/2}\max_{Q}\big(\sum_{i,j}
|\xi\delta_{h,\ell_{j}}\delta_{h,\ell_{i}}u|^{2}\big)^{1/2}
$$
$$
\leq (3\mu)^{1/2}\max_{\partial_{2}Q}\sum_{i,j}
|\xi\delta_{h,\ell_{j}}\delta_{h,\ell_{i}}u|
+ (3\mu)^{1/2}\max_{Q^{o}_{2}}\sum_{i,j}
|\xi\delta_{h,\ell_{j}}\delta_{h,\ell_{i}}u|.
$$
In light of Lemma \ref{lemma 8.24.1} the last maximum over $Q^{o}_{2}$
is less than $8d_{1}^{2}H$. Hence
$$
\theta H\leq
N\max_{\partial_{2}Q,i,j}|\xi\delta_{h,\ell_{j}}\delta_{h,\ell_{i}}u|
+(1/2)\theta H,\quad H\leq
N\max_{\partial_{2}Q,i,j}|\xi\delta_{h,\ell_{j}}\delta_{h,\ell_{i}}u| 
$$
and we can finish the proof   \eqref{10.8.1} as
 few times before.
 The theorem is proved.

\mysection{Proof of Theorem \protect\ref{theorem 11.2.1}}
                                              \label{section 11.4.1}

In the following lemma the assumption that $\sup_{\alpha}a^{\alpha}_{k}
\geq\delta$ is not used. All other assumptions 
of Theorem \ref{theorem 11.2.1} are supposed to hold. We use 
notation \eqref{5.11.1} and
the notation 
from Section \ref{section 10.7.1} with $d_{1}+1$ in place of $d_{1}$
 and 
$\delta_{i}=\delta_{h_{i},\ell_{i}}$, 
$\Delta_{i}=\Delta_{h_{i},\ell_{i}}$.
Here we take
$$
\mu=0
$$
 and show how to choose $\nu
=\nu(\lambda, \varepsilon,d_{1},K_{3})$ in Lemma \ref{lemma 11.2.1}.

\begin{lemma}
                                                 \label{lemma 11.2.1}  

In $Q|_{0}$
for $\gamma\in\Gamma(\varepsilon)$  we have
$$
( \sum_{i}\gamma_{i}\Delta_{h_{i},\ell_{i}}u)^{-}
\leq 
N^{*}e^{m^{+}(T+\tau)}(1+\max_{\bar{Q},j}
|\xi_{(-)}\delta_{h_{j},\ell_{j}} u|
$$
\begin{equation}
                                                 \label{8_27.03}
+
\max_{\bar{Q}}|\xi_{(-)}u| +
\max_{\partial_{1} Q,j}(\xi_{(-)}\Delta_{h_{j},\ell_{j}}u)^{-}),
\end{equation}
where  $N^{*}=N^{*}(\lambda,h_{0},\varepsilon, d_{1} ,K_{3})$.
\end{lemma}

Proof. As many times before, if
  $(t_{0},x_{0})\in
\partial_{1} Q$, there is nothing to prove.
Therefore, we assume that $(t_{0},x_{0})\in Q^{0}_{1}$.
We may also assume that at $(t_{0},x_{0})$
$$
 \sum_{i } 
(\xi\delta_{h_{i},\ell_{i}} u)^{2}
\leq 
\big[(\xi\sum_{i}\gamma_{0i}\Delta_{h_{i},\ell_{i}}u)^{-}\big]^{2} .
$$
 
Then
the operator $P_{\gamma_{0}0\nu}$ respects the maximum
principle for $h\nu\leq2\varepsilon$ (see the proof of Lemma
\ref{lemma 8.25.3} and recall that
$\eta\leq h$).

Then as in the proof of Lemma \ref{lemma 10.6.1}
we obtain
\begin{equation}
                                                        \label{11.4.1}
4J=4 P_{\gamma_{0}0\nu}(a _{k}\Delta_{k}
v +b _{k}\delta_{k}v)\geq  
 a _{k} \nu R^{\nu}_{k} 
 +4I_{4} 
\end{equation}
if $h\leq h_{0}/2$, where this time 
$$
I_{4}= 2v_{\gamma_{0}}^{-}\gamma_{0i}(\delta_{i}b _{k})v_{ki}
+v_{\gamma_{0}}^{-}\gamma_{0i}(\Delta_{i}b _{k})v_{k}
-
\nu v_{i}(\delta_{i}b _{k})T_{h_{i},\ell_{i}}v_{k}.
$$
Since
$$
v_{\gamma_{0}}^{-}\gamma_{0i}|(\delta_{i}b _{k})v_{ki}|
\leq N^{*}v_{\gamma_{0}}^{-}\sqrt{a _{k}}
\sum_{i}|v_{ki}|\leq (\lambda/8)(v_{\gamma_{0}}^{-})^{2}
+N^{*}a _{k}\sum_{i}v_{ki}^{2}
$$
we get from \eqref{11.4.1} that, for $\nu=
\nu^{*}(\lambda,\varepsilon,d_{1},K_{3})$,
$h \nu^{*}\leq2\varepsilon$,  and  $h\leq h_{0}/2$,
\begin{equation}
                                                     \label{11.2.3}
J\geq-(\lambda/2) (v_{\gamma_{0}}^{-})^{2}
-N^{*}\bar{W}_{1}.
\end{equation}

The rest is just a repetition of a part of the proof of Theorem
\ref{theorem 10.3.1} with obvious and big simplifications.
The lemma is proved.

There is almost nothing else to do to finish the proof
of Theorem ~\ref{theorem 11.2.1}. Indeed, \eqref{8_27.03}
with $d_{1}$ in place of $d_{1}+1$ yields the first estimate
in \eqref{11.4.4} as in the proof of Theorem \ref{theorem 10.4.1}.
After getting estimates for $|\Delta_{h,\ell_{k}}u|$
the estimate of $(\Delta_{\eta,l}u)^{-}$ follows
immediately from \eqref{8_27.03}. The theorem is proved.

\mysection{Comments on the operators having form
\protect\eqref{11.19.6}}
                                                \label{section 5.7.1}

We know (see,
for instance, \cite{DK1}) that if an operator $L$
having form \eqref{11.19.6}
 admits an approximation
with operators $S_{h}$ of the form \eqref{10.19.1}
respecting the maximum principle 
with $hB$ in place of $B$ and $\text{Span}\,B=\bR^{d}$,
 then necessarily 
$$
Lu=a _{k}\ell_{k}^{i}\ell_{k}^{j}u_{x^{i}x^{j}}
 +b  _{k}\ell_{k}^{i}u_{x^{i}}
$$
with some $a_{k},b_{k}\geq0$ and $\ell_{k}\in B\cup(-B)$.
A way to find such representations
  for $d=2$ and given $a^{ij}$ is suggested in \cite{Z}.

Next natural issue is related to the smoothness of $a_{k},b_{k}$
if we are given that $a^{ij}$ are smooth. Recall that
in Assumption \ref{assumption 10.3.1} we need $a^{\alpha}_{k}$
to be at least Lipschitz continuous. 
Of course, this problem disappears if $a^{ij}$ are constant.

It is an easy
and probably well-known fact that
if $(a^{ij})=(a^{ij}(t,x))$ is uniformly bounded and uniformly elliptic,
then one can find $d_{1}$ and $\Lambda_{0}$
for which $a_{k}$ can be chosen strictly positive
and as smooth as $a^{ij}$ are. The proof of this can be
obtained from the fact that if we are given a closed convex
polyhedron then every point in the relative 
interior can be written as
a convex combination of the extreme points with the coefficients
 $>0$ which are infinitely differentiable functions of the point.
By replacing $L$ with $L+\varepsilon^{2}\Delta$
one can approximate a possibly degenerate operator $L$ with uniformly
nondegenerate ones, so that there always exist a sequence
of operators of the form \eqref{11.19.6} approximating
$L$. 
Notice, however, that generally the set $\Lambda_{0}$ changes with
$\varepsilon$. Nevertheless, one knows how
to  estimate the difference
of solutions corresponding to $L+\varepsilon^{2}\Delta$
and $L$ (see, for instance \cite{Kr77}) and between the
solutions of the corresponding finite-difference
approximations (see, for instance, Theorem 5.6 of \cite{Kr06}
or Remark \ref{remark 11.19.3}).

Another generic example is given by the so-called diagonally dominant
matrices. For instance, take $d =2$ and assume
that $b\equiv0$, $a^{ij}$ are {\em twice\/} continuously differentiable
with respect to $x$,
 $a^{12}=a^{21}$, and
$|a^{12}|\leq s$, where
$s=a^{11}\wedge a^{22}$.

Set $\kappa=1/3$ and take an infinitely
differentiable,  even, and
 convex function $\psi(t)$ on $\bR$ such that
$\psi(y) =|y|$ for $|y|\geq\kappa$.
Introduce  
$$
g= a^{12}s^{-1} ,\quad
h=s\psi(g),
\quad
2\hat{a}^{1,\pm 2}= h\pm a^{12},\quad 
2\hat{a}^{ii}=a^{ii}-h,
$$
where $i=1,2$ and $0\cdot 0^{-1}:=0$. For other values of $i,j=\pm1,\pm2$
define $\hat{a}^{ij}$ so that 
$$
\hat{a}^{ij}=\hat{a}^{ji},\quad \hat{a}^{-i,-j}=\hat{a}^{ij},
\quad \hat{a}^{i,-i}=0
$$
 and set
$$
\ell_{j}=e_{|j|}\sign j,\quad
\ell_{ij}=\ell_{i}+\ell_{j},
$$
where $e_{1},e_{2}$ are the basis vectors. Then simple manipulations
yield
$$
4Lu=\hat{a}^{ij}\ell_{ij}^{k}\ell_{ij}^{r}
u_{x^{k}x^{r}}.
$$
 
We now show that not only $L$ admits a representation
as the sum of second order directional derivatives 
with the directions independent of $t,x$ but also
$\sqrt{\hat{a}^{ij}}$ are Lipschitz continuous in $x$.
By the way,
observe that obviously $\hat{a}^{ij}\geq0$. 

We are going to use that
nonnegative and twice continuously differentiable functions
are the squares of Lipschitz continuous functions.
In particular, $a^{ii}$  and, consequently,
$a^{11}\wedge a^{22}$ are the squares of
 Lipschitz continuous functions, 
$|a^{ii}_{x}|\leq N\sqrt{a^{ii}}$ and
$|s_{x}|\leq N\sqrt{s}$.
Furthermore,
 $a^{11}\pm a^{12}$
is   nonnegative and twice continuously differentiable. Hence,
it is the square of a Lipschitz continuous function.
In particular,  
$$
|a^{11}_{x}\pm a^{12}_{x}|\leq N\sqrt{a^{11}\pm a^{12}},
\quad | a^{12}_{x}|\leq N\sqrt{a^{11} },\quad
| a^{12}_{x}|\leq N\sqrt{s}.
$$
 and recalling that $|a^{12}|\leq s$ we find
$$
|g_{x}|\leq |a^{12}_{x}|s^{-1}
+|a^{12}|\cdot|s_{x}| s^{-2} \leq 
Ns^{-1/2},\quad|h_{x}|\leq N\sqrt{s}.
$$

Next, the function  $\phi(y) :=\psi(y)+ y$ is smooth and nonnegative.
Therefore
$$
2|\hat{a}^{12}_{x}|=|\phi'g_{x}s+\phi 
s_{x}|
\leq N\sqrt{\phi }\sqrt{s}
+N\phi \sqrt{s}
\leq N\sqrt{\phi }\sqrt{s}
=N\sqrt{\hat{a}^{12}}.
$$
Similar estimate holds for
$|\hat{a}^{1,-2}_{x}|$ and $|\hat{a}^{ij}_{x}|$ if $i\ne j$.

On the set where $|a^{12}|>\kappa s$, we have
$
2\hat{a}^{11}=a^{11}-|a^{12}|$, so that by the above
$$
|\hat{a}^{11}_{x}|=|a^{11}_{x}-a^{12}_{x}\text{sign}\,a^{12}|
\leq N\sqrt{a^{11} -a^{12} \text{sign}\,a^{12}}
=N\sqrt{\hat{a}^{11}}.
$$
Finally, on the set where $|a^{12}|<2\kappa s$
it holds that $h\leq 2\kappa s$,
$a^{11}-h\geq \kappa a^{11}$, and
 $
|\hat{a}^{11}_{x}|\leq N\sqrt{a^{11}}\leq N\sqrt{\hat{a}^{11}}$.
Similarly we get what we need for~$\hat{a}^{22}$
and the remaining $\hat{a}^{ii}$.

\end{document}